\documentclass[11pt]{amsart}
\usepackage{amssymb,amsmath,hyperref,mathabx,graphicx,float}

\newtheorem{theorem}{Theorem}

\newtheorem{conjecture}[theorem]{Conjecture}
\newtheorem{corollary}[theorem]{Corollary}
\newtheorem{definition}[theorem]{Definition}
\newtheorem{example}[theorem]{Example}

\newtheorem{fact}[theorem]{Fact}
\newtheorem{lemma}[theorem]{Lemma}
\newtheorem{problem}[theorem]{Problem}
\newtheorem{proposition}[theorem]{Proposition}
\newtheorem{question}[theorem]{Question}
\newtheorem{remark}[theorem]{Remark}

\newcommand{\bcon}{\begin{conjecture}}
\newcommand{\econ}{\end{conjecture}}
\newcommand{\bcor}{\begin{corollary}}
\newcommand{\ecor}{\end{corollary}}
\newcommand{\bdf}{\begin{definition}}
\newcommand{\edf}{\end{definition}}
\newcommand{\beq}{\begin{equation}}
\newcommand{\eeq}{\end{equation}}
\newcommand{\bexa}{\begin{example}}
\newcommand{\eexa}{\end{example}}
\newcommand{\bfac}{\begin{fact}}
\newcommand{\efac}{\end{fact}}
\newcommand{\blem}{\begin{lemma}}
\newcommand{\elem}{\end{lemma}}
\newcommand{\bprb}{\begin{problem}}
\newcommand{\eprb}{\end{problem}}
\newcommand{\bpro}{\begin{proposition}}
\newcommand{\epro}{\end{proposition}}

\newcommand{\bque}{\begin{question}}
\newcommand{\eque}{\end{question}}
\newcommand{\brem}{\begin{remark}}
\newcommand{\erem}{\end{remark}}
\newcommand{\bthm}{\begin{theorem}}
\newcommand{\ethm}{\end{theorem}}
\newcommand{\bmat}{\begin{matrix}}
\newcommand{\emat}{\end{matrix}}

\newcommand{\bpr}{\begin{proof}}
\newcommand{\epr}{\end{proof}}
  
\newcommand{\lb}{\label}

\newcommand{\comment}[1]{\,}

\newcommand{\cal}{\mathcal}
\newcommand{\p}{\partial}

\newcommand{\Z}{\mathbb Z}

\newcommand{\C}{\mathbb C}

\newcommand{\ve}{\varepsilon}

\title{Skein algebras of surfaces}

\author{J\'ozef H. Przytycki, Adam S. Sikora}

\thanks{The first author acknowledges support  
by the Simons Foundation Collaboration Grant for Mathematicians 316446.
The second author acknowledges support from U.S. National Science Foundation grants DMS 1107452, 1107263, 1107367 "RNMS: GEometric structures And Representation varieties" (the GEAR Network).}


\begin{document}

\thispagestyle{empty}

\begin{abstract} 
We show that the Kauffman bracket skein algebra of any oriented surface $F$ (possibly with marked points in its boundary) has no zero divisors and that its center is generated by knots parallel to the unmarked components of the boundary of $F$. Furthermore, we 
show that skein algebras are Noetherian and Ore. Our proofs rely on certain filtrations of skein algebras induced by pants decompositions of surfaces. We prove some basic algebraic properties of the associated graded algebras along the way.
\end{abstract}

\address{Dept. of Mathematics, George Washington University,
Washington, DC 20052, USA, and Dept. of Mathematics, Physics and Informatics, University of Gda\'nsk, Wita Stwosza 57, 80-952 Gda\'nsk, Poland}
\email{przytyck@gwu.edu}
\address{Dept. of Mathematics, 244 Math Bldg, University at Buffalo, SUNY, Buffalo, NY 14260, USA}
\email{asikora@buffalo.edu}

\maketitle


%
\section{Introduction}
%

\subsection{Skein algebras} 

Throughout the paper $(F,B)$ will denote a {\em marked surface} -- that is a connected oriented surface $F$ together with a finite subset $B$ of $\p F$ of {\em marked points}. ($\p F$ may be empty).  

An {\em $(F,B)$-link} is a framed embedding of a $1$-manifold 
$$L =S^1\cup ...\cup S^1\cup I\cup ... \cup I \hookrightarrow F\times I,$$
where $I=[0,1],$ such that the boundary points of $L$ lie in $B\times I$, and the framing of $L$ at boundary points is parallel to $\p F$. (We follow the terminology of \cite{BW1, Mu,Le2}, even though the term ``tangle'' may be more adequate here.)
We assume additionally that each arc is isotopic mod its boundary to an arc with framing parallel to $F.$ (Hence, adding a half-twist to it is not allowed.) The $(F,B)$-links are considered up to isotopy within the space of $(F,B)$-links. We denote their set, including the empty link $\emptyset,$ by $\cal L(F,B).$ 

Let $R$ be a commutative ring with identity and a fixed invertible element $A\in R.$ 
The {\em Kauffman bracket skein module $\cal{S}(F,B)$ of a marked surface $(F,B)$}  is the quotient of the free $R$-module $R \cal L(F,B)$ by the submodule generated by the Kauffman bracket skein relations (\ref{e_skein-rel}). This definition follows that of \cite{BW1} and it generalizes the skein modules of $3$-manifolds without marked points introduced by the first author and V. Turaev, cf. \cite{Bu, B-F-K, B-P,F-G,H-P, Ma,Le1,Prz-Bull, Prz, P-S,Tu, Si,S-W,Tu} and the references within. ($\cal{S}(F,\emptyset)$, denoted by $\cal{S}(F)$ is the ``standard'' Kauffman bracket skein module of $M$.)

\beq\lb{e_skein-rel}
\parbox{0.3in}{\includegraphics[height=0.3in]{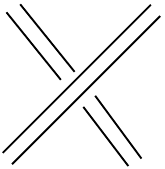}} - A
\parbox{0.3in}{\includegraphics[height=0.3in]{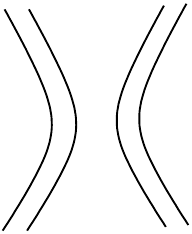}} - A^{-1}
\parbox{0.3in}{\includegraphics[height=0.3in]{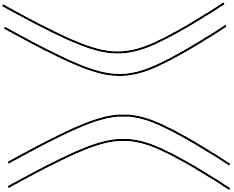}}\quad , \hspace*{.2in}
\parbox{0.3in}{\includegraphics[height=0.28in]{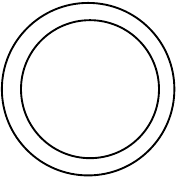}}\quad + (A^2+A^{-2})\emptyset.
\eeq

The above skein module is an algebra with the product $L_1\cdot L_2$ of two links in $\cal{S}(F,B)$ given by a union of these two links when $L_1$ is positioned in $F\times (0,1/2)$ and $L_2$ in $F\times (1/2,1).$ That multiplication operation extends onto the entire algebra $\cal S(F,B)$ with the identity $\emptyset.$ 

Skein modules and skein algebras play an important role in quantum topology.
In particular, the skein algebra of a punctured surface embeds into the quantum Teichm\"uller space of Chekhov-Fock, \cite{BW1}, cf. \cite{CF, Ka, Le2, Le3, Mu}.

Furthermore, for closed surfaces $F$, the quotient of the skein algebra of $F$ for $A=e^{2\pi i/4r}$ by the $(r-1)$st Jones-Wenzl idempotent coincides with endomorphism ring of the Reshetikhin-Turaev-Witten TQFT state space at level $r$, cf. \cite{Ro, S-tqft}. Finally, skein modules of knot complements and the skein algebra of the torus play the central role in the q-holonomicity of colored Jones polynomials, the construction of the non-commutative A-polynomial, and in the AJ-conjecture, \cite{F-G,FGL, Ga, G-L}.
 
By ``pushing'' the fibers $\{b\}\times I$ into $\p F \times \{0\}$ in the direction of the orientation of $\p F$,  $(F,B)$-links can be thought as links in $F\times I$ with endpoints in $(\p F -B)\times \{0\}$ and with framing parallel to $F$. From now on we will use this alternative definition since it makes $(F,B)$-links more conveniently represented by diagrams in $F$. By this approach the product of any two links $L_1\cdot L_2$ is given by a disjoint union of them such that $L_2$ is on top of $L_1$ and
the endpoints of $L_1$ lie before the endpoints of $L_2$ in every interval of $(\p F-B)\times \{0\}$, with respect to the orientation of $\p F.$

Let $\cal{S}'(F,B)$ be a submodule of $\cal S(F,B)$ spanned by links $L$ which contain an arc parallel to $\p F.$ (That is an arc $\alpha\subset F$ such that there is an isotopy of $L$ fixing its endpoints and moving $\alpha$ into $\alpha'\subset \p F.$ That $\alpha'$ may contain points of $B$.) Symbolically, $\cal{S}'(F,B)$ is generated by the skein relation (\ref{triv-arc}).

\beq\lb{triv-arc}
\includegraphics[height=0.3in]{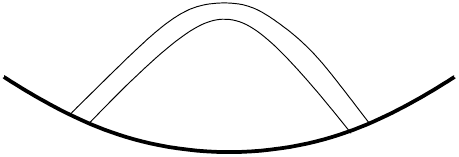}
\eeq

Note that $\cal{S}'(F,B)$ is a two-sided ideal in $\cal{S}(F,B)$. Indeed, if $L_1$ contains an arc $\alpha$ parallel to $\p F$ as above then for every $(F,B)$-link $L_2$ every resolution of crossings between $\alpha$ and $L_2$ yields a link in $\cal{S}'(F,B)$.

We call the quotient 
$${\cal RS}(F,B)={\cal S}(F,B)/{\cal S}'(F,B)$$
the {\em reduced skein algebra} of $(F,B)$. That algebra makes an appearance in \cite{Le2} already and it is a quotient of the marked skein algebra of \cite{Mu}.

For any $B\subset \p F,$ $\cal S(F)=\cal S(F,\emptyset)$ is a subalgebra of $\cal{S}(F,B)$ and since $\cal S(F)$ has trivial intersection with $\cal{S}'(F,B)$ it is also a subalgebra of $\cal{RS}(F,B)$.

For any connected component $\gamma$ of $\p F$ let $K_\gamma$ be a knot in $F\times I$ parallel to $\gamma$ (with the framing parallel to $F$). We say that a component $\gamma$ of $\p F$ is {\em marked} if it contains at least one element of $B$. (The unmarked components correspond to surface punctures in \cite{BW1,Le2,Mu}.)

Denote the center of a ring $P$ by $C(P).$ The goal of this paper is to prove the following three fundamental properties of (reduced) skein algebras.

\bthm \lb{main-cent} If $A^{4n}-1$ is not a zero divisor in $R$ for any $n>0$ then
for any marked surface $(F,B)\ne (S^1\times I,\emptyset)$, $C(\cal{RS}(F,B))$ is the $R$-algebra of polynomials in knots $K_\gamma$ for the unmarked connected boundary components $\gamma$ of $\p F.$
\ethm 

In particular, for $F$ not annulus, $C(\cal S(F))$ is a polynomial ring in $K_\gamma$ for connected boundary components $\gamma$ of $\p F$. On the other hand, 
$$C(\cal S(S^1\times I))=\cal S(S^1 \times I)=R[K],$$ 
where $K$ is the core of the annulus. A version of Theorem \ref{main-cent} was proved for some surfaces of genus zero and one (and $B=\emptyset$) in \cite{B-P}.

The condition that $A^{4n}-1$ is not a zero divisor in $R$ for $n>0$ implies that $A$ is not a root of unity.
That condition is necessary. Indeed, Bonahon and Wang prove in \cite[Lemma 10]{BW2} that for every knot $K$ in $F\times I,$ the $n$-th Chebyshev polynomial (of the first kind) of $K$ is central in $\cal S(F)$ if $A$ is a primitive $m$-th root of unity for an odd $m$. (See also \cite{Le1}.) Therefore, the center of $\cal S(F)$ is much larger than that of Theorem \ref{main-cent} in that case.

Since knots $K_\gamma$ as above are central in $\cal S(F,B)$, the epimorphism $\cal{S}(F,B)\to \cal{RS}(F,B)$ induces an epimorphism on centers of these algebras. 


The second main result of this paper is

\bthm \lb{main-div} For any marked surface $(F,B)$ and for any integral domain $R$ (i.e. a commutative ring with no non-zero zero divisors)
$\cal{RS}(F,B)$ is a domain (i.e. it has no left nor right non-zero zero divisors).  
\ethm 

These results were announced in \cite{P-S} for $B=\emptyset$.  The proofs of this paper are generalizations of those which we had in mind 18 years ago, when writing that paper.  A rough idea of the proof of Theorem \ref{main-div} was mentioned in \cite{Pr-ZD}.

In the meantime, Bonahon and Wong constructed an embedding of skein algebras, for $\p F\ne \emptyset,$ into the Chekhov-Fock quantum Teichm\"uler space, \cite{BW1,CF}, which is a quantum torus, i.e. the algebra of Laurant polynomials in non-commuting variables $x_1^{\pm 1},..., x_n^{\pm 1}$ which $q$-commute, that is $x_ix_j=q^{n_{ij}}x_jx_i,$ for some $n_{ij}\in \Z.$ Since such algebras are domains, their work implies that non-reduced skein algebras, $\cal{S}(F,B),$ are domains for $\p F\ne \emptyset$ and $R=\C[A^{\pm 1}].$ 
 (Later, Muller constructed an embedding of his skein algebra $Sk_q(F)$ for $\p F\ne \emptyset,$ into a quantum torus. Various simplifications and enhancements of Bonahon-Wang and Muller work were obtained by L\^ e, \cite{Le1, Le2, Le3}.) 

Additionally, \cite{C-M} proved that the skein algebra for $A=\pm 1$ has no non-zero nilpotents -- see further comments and consequences of this result in Subsection \ref{ss_appl_charvar}. Finally, C. Frohman and J. Kannia-Bartoszy\'nska proved the above theorem for $\p F\ne \emptyset,$ $B=\emptyset$, and $A$ root of unity, \cite{FK}.

We like to think of the above skein algebras as quantizations of algebras coming from invariant theory. One usually requires of such quantized algebras to be Noetherian, c.f. \cite[Ch. I.1.14]{BG}. As suggested to us by T. L\^e, it is indeed the case for the skein algebras:

\bthm \lb{noetherian}[Proof in Sec. \ref{ss_noetherian}] For any marked surface $(F,B)$ and for any integral domain $R$, 
$\cal{RS}(F,B)$ is Noetherian. 
\ethm

\bcor
For any marked surface $(F,B)$ and for any integral domain $R$, 
$\cal{RS}(F,B)$ is\\
(1) finitely generated\\
(2) an Ore domain, cf. \cite[Corollary 6.7]{GW}.
\ecor


\subsection{Acknowledgements} The authors would like to thank T. L\^ e and an anonymous referee for helpful comments.

\subsection{Unorientable surfaces}
For an unorientable surface $F$ one may consider an orientable twisted $[0,1]$-bundle $M$ over $F$. 
(Such a bundle is unique.)
Since there is no notion of ``up'' and ``down'' in that bundle, the previous definition of the product does not work in $M$ for generic $A$.  
However, one can still consider the product in $M$ for $A=\pm 1$ though since 
\beq \lb{e_cross-eq}
\parbox{0.3in}{\includegraphics[height=0.3in]{cross-fr.pdf}} =\parbox{0.3in}{\includegraphics[height=0.3in]{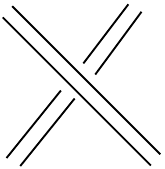}}
\eeq 
in that case and, therefore, the product of two links in $M$ can be defined as any disjoint union of them.
Let us assume that $F$ is closed now since otherwise $M$ is homeomorphic to a handlebody. 
We propose

\bcon \lb{div-unorient} 
For any nonorientable surface $F\ne \mathbb{RP}^2,\mathbb{RP}^2\hash \mathbb{RP}^2$ (Klein bottle) and for any integral domain $R$ the skein algebra $\cal S(F)$ is a domain.
\econ 

The above result does not hold for $\mathbb{RP}^2$ nor $\mathbb{RP}^2\hash \mathbb{RP}^2$.
Indeed, it is easy to see that $\cal S(\mathbb{RP}^2)\simeq R[x]/(x^2-4A)$ and $\cal S(\mathbb{RP}^2\hash \mathbb{RP}^2)\simeq R[x,y]/(x(y^2-4A)),$ if $\frac{1}{2}\in R$, cf. \cite{P-S}.

\subsection{Applications to $SL(2,\C)$-character varieties}\lb{ss_appl_charvar}
For $R=\C,$ $A=\pm 1,$  and any oriented $3$-manifold $M$ the skein module $\cal S(M)$ (which is an algebra in this case) is algebra isomorphic with the ``coordinate ring'' $\C[X(\pi_1(M))],$ of the $SL(2,\C)$-character variety,  \cite{Bu, P-S}.
More precisely, it is isomorphic with the algebra of sections of the structure sheaf of the $SL(2,\C)$-character variety, $X(\pi_1(M)),$ of $\pi_1(M)$, thought as an algebraic scheme being the geometric-invariant quotient of the $SL(2,\C)$-representation scheme of $\pi_1(M)$ by the action of $SL(2,\C)$ by conjugation. (The point being, $\C[X(\pi_1(M))]$ may a priori have nilpotents.)

\bcor For any surface $F$ the character variety $X(\pi_1(F))$ is reduced and irreducible. 
\ecor

$X(\pi_1(F))$ was proved to be reduced in \cite{C-M} and irreducible as an affine variety in \cite{R-B-C}.  A nice feature of our Theorem \ref{main-div} is that it implies the above statements as its corollaries.

As we have mentioned, $\cal S(F)$ is a domain for $R=\Z[A^{\pm 1}]$, $\p F\ne \emptyset$ by the recent work of Muller, \cite{Mu}. Furthermore, the result of \cite{C-M} together with that of \cite{R-B-C} implies no zero divisors in $\cal S(F)$ for $R=\Z[A^{\pm 1}]$ and closed surfaces, cf. Subsection \ref{ss_bases}. These results do not extend immediately however to all domains $R$ nor do they apply to $\cal{RS}(F,B).$

%
\subsection{Multi-curve bases of skein algebras}\lb{ss_bases}
%

Note that every $(F,B)$-link is represented by its diagram in $F$ -- that is a collection of (unoriented) curves, (that is circles and arcs) properly immersed in $F$ with endpoints on marked components of $F$ only, but away from $B$. 

A diagram like that without any crossings and with no contractible loops is 
a {\em multi-curve} in $(F,B)$. (In particular, the empty set is a multi-curve too.) 
We denote the set of the isotopy classes of multi-curves in $F$ (within the set of multi-curves) by $\cal{MC}(F,B)$.
By the definition, $\cal{MC}(F,\emptyset)$ contains collections of closed loops only.

By \cite{Prz,P-S, S-W},
$\cal{MC}(F,\emptyset)$ is an $R$-module basis of $\cal S(F)$. These  proofs extend verbatim to the following:

\bthm\lb{basis}
$\cal{MC}(F,B)$ is an $R$-module basis of $\cal{S}(F,B),$ for any $(F,B)$ and $R$.
\ethm

\noindent {\it Proof (sketch):} 
By resolving all crossings and eliminating contractible loops, every $(F,B)$-link can be represented as a linear combination of multi-curves. Therefore, $\cal{MC}(F,B)$ spans $\cal{S}(F,B).$ Because the Kauffman bracket skein relations are confluent in the sense of \cite{S-W} and because the multi-curves are irreducible in the sense of \cite{S-W}, they are linearly independent in $\cal{S}(F,B),$ by \cite{S-W}.
\qed

\bexa\lb{exa-D^2} 
$\cal{S}(D^2,\{*\})$ has a basis composed of crossingless matchings on finite numbers of points on the interval $\p D^2-\{*\}$, cf. Fig. \ref{crossingless}.
(Such  crossingless matchings appear in geometric representation theory, cf. for example \cite{Kh}.) 
The product on $\cal{S}(D^2,\{*\})$ defines a monoid structure on the set of crossingless matchings. It is easy to see that $\cal{S}(D^2,\{*\})$ is the free associative algebra on the set of all primitive crossingless matchings, i.e. those
which do not decompose as a product $D_1\cdot D_2$ of crossingless matchings, for $D_1, D_2\ne \emptyset.$
\eexa

We say that a multi-curve $C$ is {\em reduced} if it does not contain a boundary parallel arc. We  
denote the set of such 
multi-curves in $(F,B)$ by $\cal{RMC}(F,B).$
Since the non-reduced multi-curves form a basis of the submodule $\cal S'(F,B)\subset \cal S(F,B)$ generated by the skein relation (\ref{triv-arc}), we have

\bcor\lb{r-basis}
$\cal{RMC}(F,B)$ is an $R$-module basis of $\cal{RS}(F,B),$ for any $(F,B)$ and $R$.
\ecor

Since $\cal{RMC}(S^2,\emptyset)=\cal{RMC}(D^2,B)=\emptyset$, for any $B$, Theorems \ref{main-cent} and \ref{main-div} clearly hold for these marked surfaces.

\begin{figure}
\includegraphics[height=0.6in]{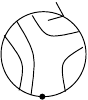}
\caption{A crossingless matching.}
\lb{crossingless}
\end{figure}

The above basis result for $\cal S(F)$ together with the results of  \cite{Bu,P-S}, \cite{R-B-C} and of \cite{C-M} provide an alternative to ours proof of Theorem \ref{main-div} for closed surfaces and $R=\Z[A^{\pm 1}]$ (but does not extend to all integral domains $R$):
\vspace*{.1in}

\noindent{\bf Proof of Theorem \ref{main-div} for $F$ closed and $R=\Z[A^{\pm 1}]$} (Suggested to us by J. Kania-Bartoszy\'nska and T. L\^e): Since $F$ is closed, $B=\emptyset.$ Consider non-zero $x,y\in \cal S(F)$. We are going to show that $x\cdot y\ne 0.$  Since closed multi-curves form a basis of $\cal S(F)$, 
$$x=\sum_C x_C\cdot C\ \text{and}\ y=\sum_C y_C\cdot C,$$
where $x_C,y_C\in \Z[A^{\pm 1}]$ and $C$ represent closed multi-curves.
By dividing $x_C$'s and $y_C$'s by a power of $A+1$ if necessary, we can assume that
at least one of $x_C$'s and at least one $y_C$'s is non-zero in
$$\cal S(F)_{\C}=\cal S(F)\otimes_{\Z[A^{\pm 1}]} \C,$$ 
where $A=-1$ in $\C.$ Since closed multi-curves are linearly independent in 
$\cal S(F)_{\C},$ the skeins $x,y$ are non-zero in $\cal S(F)_{\C}$ and, hence, by \cite{C-M}, they are also non-zero in $\cal S(F)_{\C}/\sqrt{0}.$ By \cite{Bu,P-S}, $\cal S(F,\C)/\sqrt{0}$ is the coordinate ring of the (affine) $SL(2,\C)$-character variety of $\pi_1(F)$ and, hence, $x\cdot y\ne 0$ in $\cal S(F)_{\C}$ by
\cite{R-B-C}. Consequently, $x\cdot y\ne 0$ in $\cal S(F).$ \qed

For any $R$ and $A^{\pm 1}\in R$, 
$$\cal S(F)_{R}=\cal S(F)_{\Z[A^{\pm 1}]} \otimes_{\Z[A^{\pm 1}]} R.$$ 
Therefore, it may be tempting to conclude that if $R$ is an integral domain then that tensor product has no non-zero zero divisors. However, there are known examples of integral domains whose tensor products are not integral domains.

%
\subsection{Filtered and graded skein algebras}
\label{ss_filtered_intro}
%

The basic approach to proving the results of this paper is as follows:

Consider a reduced multi-curve $M$ in $(F,B)$.
Then for every link diagram $D$ in $(F,B)$, the {\em weight} $w_M(D)$ is the minimal geometric intersection number of $D$ and $M$ in generic position. We assume that $M$ is over $D$ (as in $D\cdot M$) which means that the endpoints of $M$ lie after the points of $D$ in $\p F\setminus B.$ 
Let ${\cal F_w}$ be the subspace of ${\cal RS}(F,B)$ spanned by links with diagrams of weight $\leq w.$  Similar filtrations were considered in \cite{Mu,Le1, AF}.

\brem\lb{filtr-grad}
(1) ${\cal F_w}$ is an algebra filtration of  ${\cal RS}(F,B)$, i.e.  $\cal F_k\cdot \cal F_l\subset \cal F_{k+l}$.\\  
(2) Consequently, the above filtration yields the {\em graded reduced skein algebra}, 
$$\cal{GRS}(F,B)=\oplus_{k=0}^\infty\, \cal F_k/\cal F_{k-1},$$
where $\cal F_{-1}=\{0\}$, together with the canonical ``grading'' map 
$$gr: \cal{RS}(F,B)\to \cal{GRS}(F,B)$$
projecting elements of $\cal F_k - \cal F_{k-1}$ into $\cal F_k/\cal F_{k-1}$.  
For $x\in \cal F_k-F_{k-1},$ $y\in \cal F_l-F_{l-1},$ the multiplication $gr(x)\cdot gr(y)=(x+\cal F_{k-1})(y+\cal F_{l-1})$ in the graded algebra $\cal{GRS}(F,B)$ is given by
$xy+F_{k+l-1},$ cf. \cite[Appendix I.12]{BG}.\\
(3) Note that ${\cal RS}(F,B)$ has a basis of reduced multi-curves in $(F,B).$ Consequently,  through an identification of bases, $\cal{GRS}(F,B)$ is isomorphic to $\cal{RS}(F,B)$ as an $R$-module.
\erem

The algebra ${\cal RS}(F,B)$ can be thought as a deformed version of the graded ``simpler'' algebra, $\cal{GRS}(F,B).$

We are going to deduce Theorems \ref{main-cent}, \ref{main-div}, and \ref{noetherian} from their graded versions, cf. Theorems \ref{main-cent-gr} and \ref{main-div-gr}.\vspace*{.1in} 

%
\subsection{Proof of Theorem \ref{noetherian}}
\label{ss_noetherian}
%

As in \cite{BW1}, one can think of unmarked boundary components of $(F,B)$ as punctures.
Furthermore, for the purpose of proving that $\cal{RS}(F,B)$ is Noetherian one can freely introduce new punctures into $F$, since for a surface with new punctures $(F',B)$ one has an epimorphism $\cal{RS}(F',\emptyset)\to \cal{RS}(F,\emptyset)$.
Now the property of $\cal{RS}(F,B)$ being Noetherian can be deduced with some work from \cite{BW1}.

Here is a simple, self-contained proof of $\cal{RS}(F,B)$ being Noetherian.

By introducing a puncture if necessary one can consider an ideal triangulation of $F$, whose vertices are the punctures and the points of $B$. Let $E$ be the set of the edges of that ideal triangulation. They are either internal arcs of $F$ or segments of $\p F$ connecting marked points. Let 
$\{\cal F_w\}_{w\geq 0}$ and $\cal{GRS}(F,B)$ be the induced filtration and the graded algebra associated with $\cal{RS}(F,B)$, as in Remark \ref{filtr-grad}. By \cite[Thm. 1.6.9]{MR}, it is enough to prove that $\cal{GRS}(F,B)$ is Noetherian. By Remark \ref{filtr-grad}, $\cal{GRS}(F,B)$ has a basis composed reduced multi-curves. Since each of them is uniquely determined by its intersection numbers with the edges of the triangulation, we have an embedding $\rho: {\cal RMC}(F,B)\to \Z_{\geq 0}^E$. Its image is a set $|E|$-tuples of numbers $\{n_e\}_{e\in E}$ such that for every ideal triangle with edges $e_1,e_2,e_3$, the intersection numbers $n_{e_1},n_{e_2},n_{e_3}$ satisfy triangle inequalities and $n_{e_1}+n_{e_2}+n_{e_3}$ is even.
Hence, $\rho({\cal RMC}(F,B))$ is an additive semigroup in $\Z_{\geq 0}^E$. Although not all subsemigroups of $\Z_{\geq 0}^E$ are finitely generated, the polyhedral cones in $\Z^n$ are by \cite{Co}.
Denote the generators of $\rho({\cal RMC}(F,B))$ by $\rho(C_1),...,\rho(C_n).$ 

We say that elements $x$ and $y$ of an $R$-algebra $A$-commute if $xy=A^n yx$ for some $n\in \Z.$ Observe that all reduced multi-curves in $(F,B)$ $A$-commute in $\cal{GRS}(F,B)$. Indeed all intersections between reduced simple multi-curves $C$ and $C'$ can be pushed into the edges of the ideal triangulation. Resolutions of all such intersections result with multi-curves of lower weight, except for a single one $C''$ with $\rho(C'')=\rho(C)+\rho(C')$.
Therefore, $C\cdot C'$ equals $C''$ in $\cal{GRS}(F,B)$ up to a multiplicative power of $A$. Since for every multi-curve $M$ in $(F,B)$, $\rho(M)=\sum_{e\in E} c_e\rho(C_e),$ 
$M$ equals $\prod_{e\in E} C_e^{c_e}$ times a power of $A$. 
Since $\cal{RMC}(F,B)$ spans $\cal{GRS}(F,B)$, we conclude that $S_1,...,S_n$ generate $\cal{GRS}(F,B)$. Since $S_1,...,S_n$ $A$-commute, $\cal{GRS}(F,B)$ 
is a quantum torus. Being an iterated skew-polynomial ring, $\cal{GRS}(F,B)$ is Noetherian by \cite[Thm I.2.6]{GW}.
\qed

%
\subsection{Graded skein algebras defined by pants decompositions of surfaces}
\label{ss_graded}
%
%
Theorems \ref{main-cent} and \ref{main-div} were proved for $F$ torus in \cite{B-P}. (Note that the center of the skein algebra of the torus was proved to be trivial under slightly stricter assumptions on $R$ which are unnecessary for the proof. Theorems \ref{main-cent} and \ref{main-div} can be also deduced from the algebraic description of the skein algebra of the torus of \cite{F-G}.)

Therefore, assume from now on that $F\ne S^2,D^2, S^1\times S^1$. 
For brevity, denote the annulus $S^1\times [-1,1]$ by $\mathbb A$.
(That should not lead to the confusion with $A^{\pm 1}\in R$.)
We intend to carry proofs for $F=\mathbb A$ and $F\ne \mathbb A$ simultaneously, despite some differences between these cases. For this reason,
let $d$ be the number of boundary components of $F$, for $F\ne \mathbb A,$
and let $d=1$ for $F=\mathbb A.$ 
Denote the boundary components of $F$ by $\gamma_{1},...,\gamma_{d}$ in some arbitrary order. (For $F$ annulus $\gamma_1$ is one of the components of $\p \mathbb A$ and $\gamma_2$ is undefined.)

Now consider a maximal collection $\gamma_{d+1},...,\gamma_{d+s}$ of disjoint, non-contractible, non-parallel to each other, non-boundary-parallel, simple, closed curves in $F.$ 
If we denote the genus of $F$ by $g$, then
$$s=\begin{cases} 0 & \text{if}\ F=\mathbb A\\ 3g-3 +d & \text{otherwise,}\end{cases}$$
and if $F\ne \mathbb A$ then $F-\bigcup_{i=1}^{d+s} \gamma_i$ is a union of 
$(2g+d-2)$ disjoint open pairs of pants (the absolute value of the Euler characteristic of $F$).

$\gamma_1,...,\gamma_{d+s}$ will be fixed throughout the rest of the paper. 
Given a diagram $D$, we will denote the minimal geometric intersection number between $D$ and $\gamma_i$  by $n_i(D)$ for $i=1,...,d+s.$
From now on, the {\em weight} of $D$ will refer to the weight with respect to $M=\bigcup_{i=1}^{s+d} \gamma_i$, i.e. $\sum_{i=1}^{d+s} n_i(D).$
We are going to consider the filtration $\{\cal F_k\}_{k\geq 0}$ and the graded
skein algebras $\cal{GRS}(F)$ introduced in Sec. \ref{ss_filtered_intro} (with respect to
$M=\bigcup_{i=1}^{s+d} \gamma_i$).  
We are going to deduce Theorems \ref{main-cent} and \ref{main-div} from their graded versions:

\bthm \lb{main-cent-gr} 
If $A^{4n}-1$ is not a zero divisor in $R$ for any $n>0$ then for any marked surface $(F,B)\ne (\mathbb A,\emptyset)$, $C(\cal{GRS}(F,B))$ is the polynomial ring (with coefficients in $R$) in variables $K_\gamma$ for the unmarked connected components $\gamma$ of $\p F.$
\ethm

\bthm \lb{main-div-gr} 
If $R$ is an integral domain then for every marked surface $(F,B)$ the graded reduced skein algebra $\cal{GRS}(F,B)$ is a domain.
\ethm

The above theorem immediately implies Theorem \ref{main-cent} and
 \ref{main-div}. Let us start with the latter one: Let $x\in \cal F_k-\cal F_{k-1},$
 $y\in \cal F_l-\cal F_{l-1}$, for some $k,l\geq 0.$ Then $gr(x)\cdot gr(y)=x\cdot y+\cal F_{k+l-1}$ is non-zero in $\cal{GRS}(F,B)$ implying that $x\cdot y\not\in \cal F_{k+l-1}$ and, hence, non-zero.

For the first one, notice that if $c\in C(\cal{RS}(F,B))$, $c\ne 0$ then $gr(c)$ is a non-zero element of the center of $\cal{GRS}(F,B)$. By Theorem \ref{main-cent-gr}, $gr(c)$ is a polynomial in $K_\gamma$ for the unmarked connected components $\gamma$ of $\p F.$
Consequently, $gr(c)\in \cal F_0$ and, hence, $c=gr(c).$ This implies that $C(\cal{RS}(F,B))$ is generated by $K_\gamma$'s as above. They are algebraically independent, because monomials in these variables are basis elements of $\cal{RS}(F,B)$ by Theorem \ref{basis}.

%
\section{Classification of multi-curves in surfaces}
\lb{s_class}
%

Let $\mathbb A_1,...,\mathbb A_{d+s}$ be mutually disjoint regular neighborhoods of $\gamma_1,...,\gamma_{d+s}$.
For $F=\mathbb A$ (annulus), take $\mathbb A_1=\mathbb A$.

\begin{figure}
\includegraphics[height=1in]{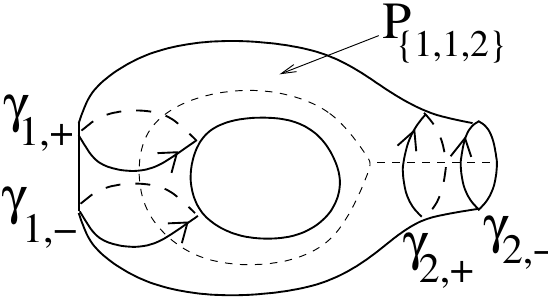}
\caption{Decomposition of once-punctured torus into annuli and a pair of pants. Base intervals in annuli and the $Y$ graph in the pair of pants are depicted by the dash lines in the back.}
\lb{pants-dec}
\end{figure}

Fix an orientation and a basepoint $*$ in $S^1.$  Identify each $\mathbb A_i$ with $\mathbb A=S^1\times [-1,1]$ through an orientation preserving homeomorphism. 
Let $\gamma_{i,+},$ $\gamma_{i,-}$ be the components of $\p \mathbb A_i$ corresponding 
to $S^1\times \{+1\}$ and $S^1\times \{-1\}$ respectively. We orient $\gamma_{i,\pm}$ in the direction of $S^1$.
Note that by the above conventions $\gamma_{i,-}=\gamma_{i}\subset \p F$ for $i=1,...,d.$

We call the image of $\{*\}\times [-1,1]\subset \mathbb A$ under the identification $\mathbb A_i\simeq \mathbb A$ the {\em base interval of $\mathbb A_i$}, for $i=1,...,d+s.$ We will assume that if $\mathbb A_i$ is adjacent to a marked boundary component of $F$ then the end of the base interval in $\mathbb A_i$ lying in $\p F$ belongs to $B$. 

Since $F\ne S^2,D^2,S^1\times S^1,$ the above annuli decompose $F$ into pairs of pants.  (For $F=\mathbb A$ the set of pairs of pants is empty.) Every pair of pants in $F$ is bounded by $\gamma_{i,\ve_1},\gamma_{j,\ve_2},\gamma_{k,\ve_3},$ for some $1\leq i,j,k\leq d+s$, $\ve_1,\ve_2,\ve_3\in\{\pm \}.$ We denote it by $\mathbb P_{\{i,j,k,\ve_1,\ve_2,\ve_3\}}$ or by $\mathbb P_{\{i,j,k\}}$ for simplicity.
 The index $\{i,j,k\}$ is a multi-set since two of the indices $i,j,k$ may coincide. Denote the set of all such multi-sets $\{i,j,k\}$ 
by $\cal P.$\footnote{Note that $\{i,j,k\}$ determines $\mathbb P_{\{i,j,k\}}$ uniquely, with the only exception being when two pairs of pants are glued together into a closed surface $F$ of genus $2$. Despite this, we will use notation $\mathbb P_{\{i,j,k\}}$ in our proofs for simplicity. All arguments work for the above special case as well, by simply including $\ve_1,\ve_2,\ve_3$ in the subscripts.}
 
Let $\mathbb P$ be a model pair of pants. We assume that $\mathbb P$ is oriented and that it contains a {\em base graph} $Y$  composed of three edges with one endpoint on each of the components of $\p \mathbb P$ and the other endpoints meeting a trivalent vertex.  (Hence, the notation, $Y$.) We identify each $\mathbb P_{\{i,j,k\}}$ with $\mathbb P$ through some homeomorphism such that the endpoints $Y\cap \mathbb P_{\{i,j,k\}}$
coincide with the endpoints of the base intervals of the three surrounding annuli $\mathbb A_i,\mathbb A_j,\mathbb A_k,$ cf. Fig \ref{pants-dec}.

Our proofs rely on the classification of reduced multi-curves on surfaces based on the above parametrized pants decomposition, due to Dehn-Thurston, \cite{De, F-L-P, P-H}, which was extended to multi-curves containing arcs in \cite{P-P}. (That ``parametrization'' is determined by the graph in $F$ composed of the base intervals and Y-graphs defined above.)
According to that classification, given a parametrized pants decomposition of $F$ (as above), all reduced multi-curves in $F$ can be positioned in a unique canonical form composed of {\em canonical diagrams} in each annulus and in each pair of pants. These canonical diagrams are defined as follows.

For any $(n,t) \in \Z_{> 0}\times \Z$, let $T(n,t)$ be a $1$-manifold properly embedded in $\mathbb A$ composed of $n$ arcs connecting $S^1\times \{1\}$ with $S^1\times \{-1\}$, away from $\{*\}\times \{\pm 1\}$, 
such that the embedding realizes minimal geometric intersection number $|t|$ between $T(n,t)$ and $\{*\}\times [-1,1]$ and that each intersection of $T(n,t)$ with $\{*\}\times [-1,1]$ is in the direction of $S^1$ for $t>0$ as the ``height'' coordinate in $[-1,1]$ increases, or in the opposite direction for $t<0$, cf. for example Fig. \ref{tangle-anulus}. (For consistency, in all drawings of surfaces we assume that the most upfront piece is oriented counterclockwise.) There is an alternative way of distinguishing between positive and negative twists in annuli -- a positive twist is obtained by the $A$ smoothing of the crossing in the multiplication of the base interval by the meridian.\footnote{That means the meridian is on top. The remark about the $A$-smoothing implies that orientating $\gamma_{i,\ve}$ curves is not essential for the proofs. It will be sometimes useful however for bookkeeping.}
\begin{figure}
\includegraphics[height=1in]{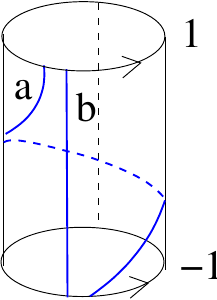}
\caption{Canonical tangle diagram $T(a+b,a)$ in $\mathbb A$ for $a,b\geq 0$ which denote multiplicities of parallel arcs.}
\lb{tangle-anulus}
\end{figure}

Furthermore, for $t\geq 0$, let $T(0,t)$ be a disjoint collection of $t$ simple closed loops in $\mathbb A$ parallel to $\p \mathbb A$.
The diagrams $T(n,t)\subset \mathbb A$ for $(n,t)\in \Z_{> 0}\times \Z\cup \{0\}\times \Z_{\geq 0}$ are the {\em canonical diagrams} in $\mathbb A$. 

We define ``canonical'' diagrams $T(n_0,n_1,n_2)$ in $\mathbb P$ as follows: 
Denote the boundary components of $\mathbb P$ by $\p_0 \mathbb P,$ $\p_1 \mathbb P,$ and $\p_2 \mathbb P$ as in Fig. \ref{pants-t}(a).  Then for any $n_0,n_1,n_2\in \Z_{\geq 0}$ with $n_0+n_1+n_2$ even, $T(n_0,n_1,n_2)\subset \mathbb P$ is defined as follows: 

(a) If $n_0,n_1,n_2$ satisfy triangle inequalities,
\beq\lb{e_triang-ineq}
n_0\leq n_1+n_2,\quad n_1\leq n_0+n_2,\quad n_2\leq n_0+n_1,
\eeq
then $T(n_0,n_1,n_2)$ is composed of 
$n_i'=\frac{n_{i+1}+n_{i+2}-n_i}{2}$ parallel arcs connecting $\p_{i+1} \mathbb P$ with $\p_{i+2} \mathbb P,$
where $i=0,1,2$ and all indices are considered mod $3.$
These arcs are disjoint and none of them intersects $Y,$ cf. Pic \ref{pants-t}(a). 
\begin{figure}
\includegraphics[height=1.3in]{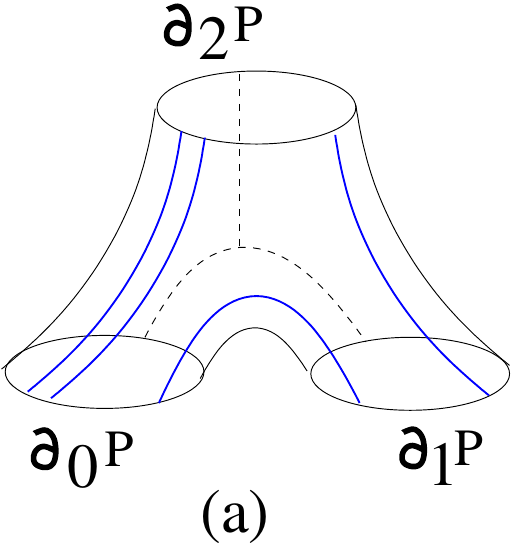}\quad
\includegraphics[height=1.3in]{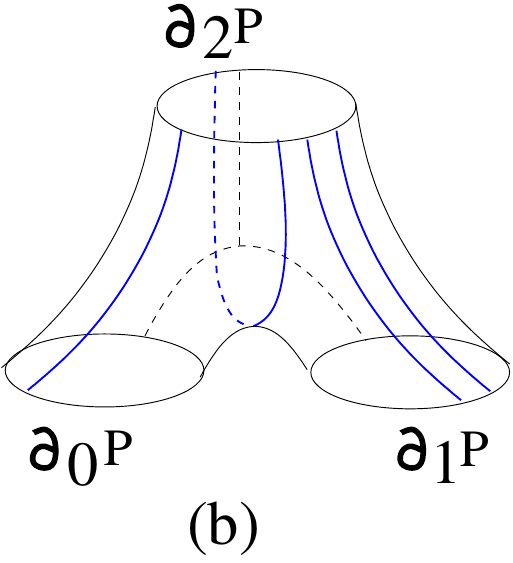}\quad
\includegraphics[height=1.3in]{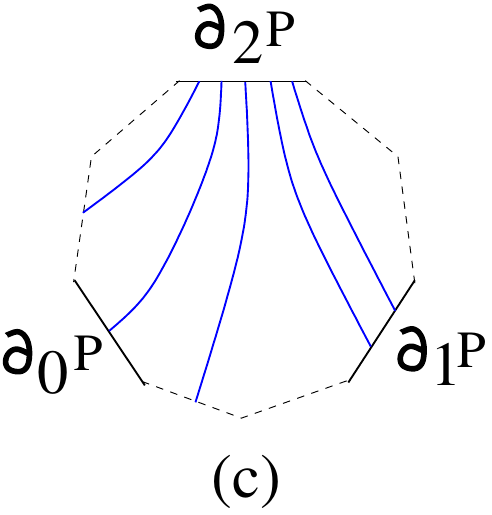}
\caption{Canonical tangle diagrams for the pair of pants $\mathbb P$. The base graph $Y$ is drawn by dashed lines in the back. (a)  $T(3,2,3)$ (b)  $T(1,2,5)$ (c) $T(1,2,5)$ drawn on a pair of pants cut along the graph 
$Y.$}
\lb{pants-t}
\end{figure}

(b) If $n_0,n_1,n_2$ fail to satisfy triangle inequalities, then $n_i>n_{i+1}+n_{i+2}$, for some $i\in \{0,1,2\}$
and we define $T(n_0,n_1,n_2)$ to be 
composed of $n_j$ parallel arcs connecting $\p_i \mathbb P$ and $\p_j \mathbb P,$ for $j=i+1,i+2$ mod $3$, 
and of 
$$n_i'=(n_i-n_{i+1}-n_{i+2})/2$$ 
parallel arcs $\alpha$ connecting $\p_i \mathbb P$ with itself and going ``around'' $\p_{i+1} \mathbb P.$ 
That is if you cut $\mathbb P$ along $Y$ and place it on a plane so that 
it has the counter-clockwise orientation, then $\alpha$ goes around $\p_{i+1} \mathbb P$, 
cf. Fig. \ref{pants-t}(b,c).
As before all of the above arcs are disjoint from each other.

Note that by this definition, $T(n_0,n_1,n_2)$ has precisely $n_i$ endpoints on $\p_i \mathbb P$, $i=1,2,3.$

Let $B'\subset \p F$ be the set of endpoints of base intervals in $\mathbb A_1,...,\mathbb A_d$ lying in the marked components of $\p F$. By our earlier assumption about endpoints of base intervals, $B'\subset B.$
We say that a reduced multi-curve in $\cal{RMC}(F,B')$ is in a canonical position iff its intersection with each annulus $\mathbb A_i$ and with each pairs of pants in $\cal P$ is a canonical diagram.
By Dehn-Thurston-Papadopoulos-Penner classification theorem, any reduced multi-curve is isotopic to a unique multi-curve in a canonical position, \cite{P-P}. 

Recall that for any reduced multi-curve $C\in \cal{RMC}(F,B')$, $n_1(C),...,n_{d+s}(C)$ are the geometric intersection numbers of $C$ with $\gamma_1,...,\gamma_{d+s}$. Additionally, let $t_1(C),....,t_{d+s}(C)$ be defined as follows: If $C'$ is the unique canonical multi-curve isotopic to $C$ then $C'\cap \mathbb A_i=T(n_i(C),t_i(C))$ for $i=1,...,d+s$. 

These are {\em Dehn-Thurston numbers (or coordinates)} of $C$. The ``$n$''-s are called the ``intersection''  numbers and the ``$t$''-s are the ``twisting'' numbers.
Hence, we arrive at

\bthm\lb{D-T-thm}
For a surface $(F,B')$ with a parametrized pants decomposition, reduced multi-curves $C$ in $(F,B')$ are classified by intersection numbers 
$n_1(C),...,n_{d+s}(C)\in \Z_{\geq 0},$ 
and by the ``twisting'' numbers, $t_1(C),...,t_{d+s}(C)\in \Z.$
\ethm

Obviously not all $2(d+s)$-tuples of intersection and twist numbers are realized for given $(F,B').$ They are intersection and twist numbers of some multi-curve in $(F,B')$ iff (a) for any pair of pants the surrounding intersection numbers adds up to an even number, (b) if $n_i=0$ then $t_i\geq 0$, (c) $n_{i}=0$ for every unmarked boundary component $\p_i F.$
 
The classification of multi-curves in all marked surfaces $(F,B)$ for arbitrary $B$
is exactly as above, except that one needs to take into account different possible positions of multi-curve  endpoints with respect to points of $B$. Specifically, we have:

\bcor\lb{D-T-B-cor}
Reduced multi-curves $C$ in $(F,B)$ are classified by intersection numbers 
$n_b(C)$ for $b\in B$, which count the intersection numbers of $C$ with intervals of $\p F$ with initial point $b$ (with respect to the orientation of $\p F$), the intersection numbers
$n_{d+1}(C),...,n_{d+s}(C),$ 
and the ``twisting'' numbers, $t_1(C),...,t_{d+s}(C)\in \Z.$ 
\ecor

%
\section{Multiplications in annuli}
%

Recall that $\mathbb A=S^1\times [-1,1]$ and that $S^1$ has a fixed base point $*.$ Let $\cal S_n(\mathbb A)\subset \cal S(\mathbb A,\{*\}\times \{\pm 1\})$ be composed of link diagrams with precisely $n$ endpoints in $S^1\times \{1\}$ and $n$ endpoints in 
$S^1\times \{-1\}$. Recall that $\cal S'(F,B)$ denotes a submodule of $\cal S(F,B)$ spanned by link diagrams containing a boundary parallel arc. By analogy, we denote the submodule of $\cal S_{n}(\mathbb A)$ spanned by such link diagrams by
$\cal S_{n}'(\mathbb A)$ and we define $\cal{RS}_n(\mathbb A)$ to be $\cal S_{n}(\mathbb A)/\cal S_{n}'(\mathbb A)$.

$\cal{RS}(\mathbb A,\{*\}\times \{\pm 1\})$ has a basis composed of reduced multi-curves in $(\mathbb A,\{*\}\times \{\pm 1\})$ with no caps nor cups. These multi-curves have equal numbers of top and bottom endpoints and, consequently,
$$\cal{RS}(\mathbb A,\{*\}\times \{\pm 1\})=\oplus_{n\geq 0}\, \cal{RS}_n(\mathbb A).$$
In fact, this is an algebra grading on $\cal{RS}(\mathbb A,\{*\}\times \{\pm 1\})$.
(Note however that $\cal{S}(\mathbb A,\{*\}\times \{\pm 1\})$ is much larger than
$\oplus_{n\geq 0}\, \cal{S}_n(\mathbb A).$)

For the purpose of proving formulas for products $T(n_1,t_1)\cdot T(n_2,t_2)$ in skein algebras of the annulus, it is useful to consider a concatenation operation
$$\cal S_n(\mathbb A)\times \cal S_n(\mathbb A)\xrightarrow{\circ} 
\cal S_n(\mathbb A)$$
obtained by identifying $S^1\times \{-1\}$ of the first annulus with $S^1\times \{1\}$ of the second one and by gluing the corresponding endpoints of link diagrams.
(We denote this concatenation multiplication by $\circ$ to distinguish it from the skein algebra multiplication.)
Note that
\beq\lb{e_annulus-tang-comp}
T(n,t)\circ T(n,t')=T(n,t+t').
\eeq 

One may call $(\cal{S}_{n}(\mathbb A),\circ)$ the ``Temperley-Lieb algebra of the annulus''.
 
\blem  $\cal S_{n}'(\mathbb A)$ is a two-sided ideal in $\cal S_{n}(\mathbb A)$ with respect to $\circ$ and, consequently, $(\cal{RS}_{n}(\mathbb A),\circ)$ is an algebra. 
\elem

\bpr
$\cal S_{n}'(\mathbb A)$ is spanned by diagrams with a cap or a cup.
$\cal S_{n}(\mathbb A)$ is spanned by multi-curves.
A concatenation of a diagram with a cap or cup with a multi-curve contains a cap or a cup and, hence, it belongs to $\cal S_{n}'(\mathbb A)$ as well. 
\epr

\brem \lb{annulus-circ}
(1) $\cal{RS}_{n}(\mathbb A)$ is a free $R$-module with a basis given by $T(n,t)$ for $t\in \Z$.\\
(2) The products $\cdot$ and $\circ$ are ``distributive'' with respect to each other in the following sense: if $x_i,y_i\in \cal{RS}_{n_i}(\mathbb A)$ for $i=1,2,$ then
$$(x_1\circ y_1)\cdot (x_2\circ y_2)=(x_1\cdot x_2)\circ (y_1\cdot y_2)$$
in $\cal{RS}_{n_1+n_2}(\mathbb A).$
\erem 

\blem \lb{lem-elem-tang-prod} 
If $n_1,n_2>0$, then
$$T(n_1,t_1)\cdot T(n_2,t_2)=A^{t_2 n_1-t_1 n_2}\cdot T(n_1+n_2,t_1+t_2),$$
in $\cal{RS}_{n_1+n_2}(\mathbb A)$.
\elem 

The statement holds for $n_1=n_2=0$ as well, but not when precisely one of $n_1, n_2$ is zero, cf. equations (\ref{e-annulus-l}) and (\ref{e-annulus-r}).

\bpr 
By Remark \ref{annulus-circ}(2),
$$T(n_1,t_1)\cdot T(n_2,t_2)=(T(n_1,t_1)\circ T(n_1,0))\cdot (T(n_2,0)\circ T(n_2,t_2))=$$ 
$$(T(n_1,t_1)\cdot T(n_2,0))\circ (T(n_1,0)\cdot T(n_2,t_2)).$$
Therefore, since twists are additive under composition $\circ,$ it is enough to prove that
\beq \lb{e_TT}
T(n_1,t_1)\cdot T(n_2,0)=A^{-t_1\cdot n_2}\cdot T(n_1+n_2,t_1)
\eeq 
and
$$T(n_1,0)\cdot T(n_2,t_2)=A^{t_2\cdot n_1}\cdot T(n_1+n_2,t_2).$$

We prove the first identity, (\ref{e_TT}), only, since the proof of the second one is practically identical. Assume that $t_1\ne 0$, since otherwise (\ref{e_TT}) is obvious.
Note that $T(n_1,\ve)\cdot T(n_2,0)$, for $\ve =t_1/|t_1|\in \{\pm 1\},$ has an $n_2$ crossing diagram and all resolutions of these crossings produce diagrams in $\cal S_{n_1+n_2}'(\mathbb A)$ with the exception of the one
obtained by $A^{-\ve}$ smoothings only (cf. Fig. \ref{fig-resolution}). That resolution yields diagram $T(n_1+n_2,\ve).$
Hence, 
$$T(n_1,\ve)\cdot T(n_2,0)= A^{-\ve\cdot n_2}\cdot  T(n_1+n_2,\ve).$$
By composing this identity $|t_1|$ times with itself (by $\circ$)
we obtain the desired identity (\ref{e_TT}).
\epr


\begin{figure}
\parbox{0.7in}{\includegraphics[height=0.7in]{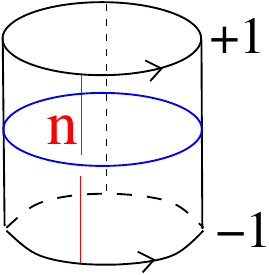}} $=
A^n$ \parbox{0.7in}{\includegraphics[height=0.7in]{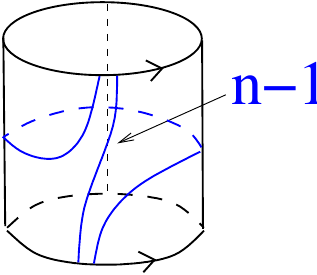}} $+ A^{-n}$
\parbox{0.7in}{\includegraphics[height=0.7in]{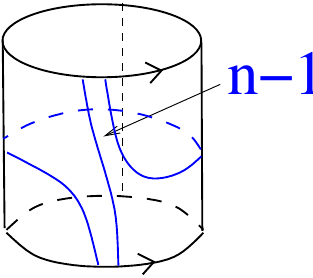}}
\caption{$T(0,1)\cdot T(n,0)$ in $\cal{RS}_n(\mathbb A)$.}
\lb{annulus-tg}
\end{figure} 

\blem \lb{lemma-annulus-tg} 
$$T(n,t)\cdot T(0,1)=A^{n}\cdot T(n,t+1)+A^{-n}\cdot T(n,t-1)\quad  \text{and}$$ 
$$T(0,1)\cdot T(n,t)=A^{-n}\cdot T(n,t+1)+A^{n}\cdot T(n,t-1)$$ 
in $\cal{RS}_{n}(\mathbb A)$,
cf. Fig. \ref{annulus-tg}. 
\elem 

\begin{figure}
Resolutions of \parbox{0.4in}{\includegraphics[height=0.4in]{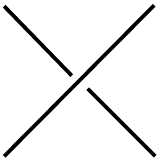}} :\quad
$A$-type:\quad \parbox{0.4in}{\includegraphics[height=0.4in]{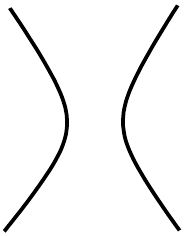}}\ ,\quad
$A^{-1}$-type:\quad \parbox{0.4in}{\includegraphics[height=0.4in]{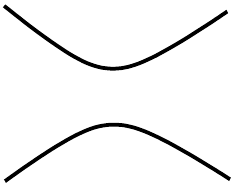}}\ .
\caption{Resolutions of a crossing}
\lb{fig-resolution}
\end{figure}

\bpr 
Since $T(n,t)=T(n,0)\circ T(n,t)$ and $T(0,1)= T(0,1)\circ \emptyset$, Remark \ref{annulus-circ}(2) implies that it is enough to prove the identities of the statement for $t=0.$

All resolutions of the $n$ crossings of $T(n,0)\cdot T(0,1)$ and of $T(0,1)\cdot T(n,0)$ result in cups or caps and, hence, diagrams in $\cal S_{n}'(\mathbb A)$, except for the one with all resolutions of type $A$ and the one with all resolutions of type $A^{-1}$ which yield the equations of the statement, cf. Fig. \ref{annulus-tg}.
\epr

Since $T(0,t')$ is $T(0,1)$ $\circ$-composed with itself $t'$ times, an induction based on the above lemma yields
\beq\lb{e-annulus-l}
T(n,t)\cdot T(0,t')= \sum_{j=0}^{t'} {t' \choose j} A^{(t'-2j)n}\, T(n,t+t'-2j)
\eeq
and
\beq\lb{e-annulus-r}
T(0,t')\cdot T(n,t)= \sum_{j=0}^{t'} {t' \choose j} A^{-(t'-2j)n}\, T(n,t+t'-2j)
\eeq
in $\cal{RS}_{n}(\mathbb A)$.

%
\section{Graded skein algebras}
%

By analogy with the subspaces $\cal{S}_n(\mathbb A)$ of $\cal{S}(\mathbb A)$ consider subspaces $\cal{S}_{n_0,n_1,n_2}(\mathbb P)$ of $\cal{S}(\mathbb P,\p \mathbb P\cap Y)$ composed of linear combinations of $(\mathbb P,\p \mathbb P\cap Y)$-link diagrams with $n_i$ endponts on $\p_i \mathbb P$ boundary component of $\mathbb P$ for $i=0,1,2.$
 We define $\cal{RS}_{n_0,n_1,n_2}(\mathbb P)$ in an analogous way, by killing all diagrams with boundary parallel components.
Clearly, 
$$\cal{S}(\mathbb P,\p \mathbb P\cap Y)=\oplus_{n_0,n_1,n_2\geq 0}\, \cal{S}_{n_0,n_1,n_2}(\mathbb P)$$ 
and an analogous decomposition holds for $\cal{RS}(\mathbb P,\p \mathbb P\cap Y)$ as well.

For every sequence of numbers $n_1,...,n_{d+s}\geq 0$ our parametrized pants decomposition 
$$F=\bigcup_{i=1}^{d+s} \mathbb A_i\cup \bigcup_{\{i,j,k\}\in \cal P} \mathbb P_{\{i,j,k\}}$$
defines a concatenation map
$$\Phi: \bigtimes_{i=1}^{d+s} \cal{S}_{n_i}(\mathbb A) \times \bigtimes_{\{i,j,k\}\in \cal P} 
\cal{S}_{n_i,n_j,n_k}(\mathbb P)\to \cal{S}(F,B'),$$
where as before $B'$ consists of endpoints base intervals in the marked components of $\p F$. Here, $n_1,...,n_{d+s}$ are any non-negative numbers (with $n_i=0$ if $i\leq d$ and $\p_i F$ is unmarked).

\blem\lb{l_conc}
The above concatenation map descends to
$$\Phi': \bigtimes_{i=1}^{d+s} \cal{RS}_{n_i}(\mathbb A) \times \bigtimes_{\{i,j,k\}\in \cal P} 
\cal{RS}_{n_i,n_j,n_k}(\mathbb P)\to \cal F_w/\cal F_{w-1}\subset \cal{GRS}(F,B'),$$
for $w=\sum_{i=1}^{d+s} n_i.$
\elem

\bpr 
Consider any link diagrams in $\cal S_{n_1}(\mathbb A),...,\cal S_{n_{d+s}}(\mathbb A)$ and in
$\cal S_{n_i,n_j,n_k}(\mathbb P)$, for $\{i,j,k\}\in \cal P.$ If one of them belongs to 
$\cal S_{n_1}'(\mathbb A),...,\cal S_{n_{d+s}}'(\mathbb A)$ or to some $\cal S_{n_i,n_j,n_k}'(\mathbb P)$ then the link obtained through contractions $\Phi$ either (a) can be deformed by pushing a boundary parallel arc (``cap'') through one of the gamma curves, which will result in lowering its intersection number with that curve or (b) it contains an arc parallel to $\p F.$
\epr

\brem\lb{r-gen-con} 
Lemma \ref{l_conc} generalizes easily to marked surfaces $(F,B)$ for any $B$. In that case, for any component $\p_i F$ of $\p F$ with $\p_i F\cap B=\{b_1,...,b_k\}$ and $k>1,$ one needs to consider $\cal{RS}_{n_i}(\mathbb A, \{*\}\times \{ 1\}\cup \{b_1,...,b_k\}\times \{-1\}\})$ (instead of $\cal{RS}_{n_i}(\mathbb A)$), 
where $\{- 1\}\times S^1$ is identified with $\p_i F$.
\erem

More generally, consider marked surfaces surfaces $(F_1,B_1),...,(F_k,B_k)$ with some parametrized pants decompositions. Gluing them along some boundary components with compatible numbers of base points results in a marked surface $(F,B)$ together with a corresponding concatenation map
$$\cal{GRS}(F_1,B_1)\times ... \times \cal{GRS}(F_k,B_k)\to \cal{GRS}(F,B).$$

%
\section{Proof of Theorem \ref{main-cent-gr}}
\lb{s_proof-cent}
%


The knots $K_\gamma$ for the unmarked components $\gamma$ of $\p F$ are obviously central in 
$\cal{GRS}(F,B)$. Since $\cal{GRS}(F,B)$ has a basis of reduced multi-curves, monomials in knots $K_\gamma$ are linearly independent. Hence, $K_\gamma$'s are algebraically independent. Consequently, to complete the proof of the statement it is enough to show that $K_\gamma$'s generate $C(\cal{GRS}(F,B)).$ 

Recall that by Corollary \ref{D-T-B-cor} the reduced multi-curves in $(F,B)$ are classified by
\begin{itemize}
\item intersection coordinates $n_b,$ for $b\in B,$ which we will think of as a function $\nu: B\to \Z_{\geq 0}$, $\nu(b)=n_b,$
\item intersection coordinates $n_{d+1},...,n_{d+s}\geq 0,$ and
\item the twist coordinates $t_1,...,t_{d+s}.$
\end{itemize}

Let us assume that $\gamma_1,...,\gamma_u$ are the unmarked components of $\p F,$ while 
$\gamma_{u+1},...,\gamma_{d}$ are the marked ones.  For $F$ annulus, for which we assume $B\ne \emptyset,$ let $u=1$ if one of the boundary curves of $F$ is unmarked and $u=0$ if both are marked. Then $t_1,...,t_u\geq 0$ for multi-curves in $(F,B).$

Let $\cal H_{\nu, n_{d+1},...,n_{d+s}}$ be a free $R$-module with a basis given by multi-curves in $\cal{RMC}(F,B)$ with intersection coordinates $\nu:  B\to \Z_{\geq 0}$, $n_{d+1},...,n_{d+s}$. We say that $\nu, n_{d+1},...,n_{d+s}$ are {\em realizable} iff $\cal H_{\nu, n_{d+1},...,n_{d+s}}\ne \{0\}$.
We call such submodules of $\cal{GRS}(F,B)$ {\em $\cal H$-spaces}.
Each $\cal F_k/\cal F_{k-1}$ decomposes into a direct sum of $\cal H$-spaces for $\sum_{b\in B} \nu(b)+n_{d+1}+...+n_{d+s}=k$. Furthermore, this decomposition makes $\cal{GRS}(F,B)$ into a multi-graded algebra. Therefore, any central element of $\cal{GRS}(F,B)$ 
is a sum of central elements in spaces $\cal H_{\nu, n_{d+1},...,n_{d+s}}$ and, consequently, 
it is enough to prove that every $c\in \cal H_{\nu, n_{d+1},...,n_{d+s}}$ central in $\cal{GRS}(F,B)$ is a polynomial in $K_\gamma$'s.

For any realizable $\nu, n_{d+1},...,n_{d+s}$ identify $\cal H_{\nu, n_{d+1},...,n_{d+s}}$ with\\ $R[x_1,...,x_u, x_{u+1}^{\pm 1},..., x_{d+s}^{\pm 1}],$ by assigning to every reduced multi-curve $C$ the monomial $x_1^{t_1(C)},...
,x_{d+s}^{t_{d+s}(C)}.$
This identification is very convenient for expressing certain products in $\cal{GRS}(F,B)$. Consider the following example:

\brem\lb{rem_mult_gamma}
By Remark \ref{r-gen-con}, 
the left and right multiplications of a reduced multi-curve $C$ in $\cal{GRS}(F,B)$ by $\gamma_i$ affects the $i$-th twisting number of $C$ only.
Specifically, by Lemma \ref{lemma-annulus-tg},  $C\cdot \gamma_i$ and  $\gamma_i\cdot C$
are represented by 
$$x_1^{t_1(C)}\cdot ... \cdot x_{d+s}^{t_{d+s}(C)}\cdot (A^{n_i} x_i+A^{-n_i}x_i^{-1})$$ 
and 
$$x_1^{t_1(C)}\cdot ... \cdot x_{d+s}^{t_{d+s}(C)}\cdot (A^{-n_i} x_i+A^{n_i}x_i^{-1})$$
in $\cal H_{n_1,...,n_{d+s}}$.
\erem

Assume that $A^{4n}-1$ is not a zero divisor in $R$ for any $n\ne 0.$
Consider a non-zero $c\in \cal H_{\nu, n_{d+1},...,n_{d+s}}\cap C(\cal{GRS}(F,B))$ as above for some $\nu, n_{d+1},..., n_{d+s}.$

\blem
$\nu=n_{d+1} = ... = n_{d+s}=0.$
\elem

\bpr Assume $n_i\ne 0$ for some $i\in B\cup \{d+1,...,d+s\}$.
Since $c$ is an $R$-linear combination of reduced multi-curves in $\cal H_{\nu, n_{d+1},...,n_{d+s}}$, by Remark \ref{rem_mult_gamma},
$$c\cdot \gamma_i-\gamma_i\cdot c=c\cdot (A^{n_i}-A^{-n_i}) (x_i-x_i^{-1})$$
in $R[x_1,...,x_u, x_{u+1}^{\pm 1},..., x_{d+s}^{\pm 1}].$
(As before, $n_i$ for $i\leq d$ are sums of $n_b$ over all points $b$ on the $i$-th boundary component.)
This expression is non zero since 
$$c\cdot (A^{n_i}-A^{-n_i})=c\cdot A^{-n_i}(A^{2n_i}-1).$$ 
Hence, $c\cdot \gamma_i\ne \gamma_i\cdot c$ -- yielding a contradiction.
\epr

Therefore, $c$ is a polynomial expression in curves $\gamma_1,...,\gamma_{d+s}.$ 
We need to prove that $c$ involves variables $\gamma_1,...,\gamma_u$ only. 
Assume that is not the case and $c$ involves a monomial 
$$m=r\cdot \gamma_1^{t_1}\cdot ...\cdot \gamma_{d+s}^{t_{d+s}},$$
with $t_i>0$ for some $i>u.$

Let $\Omega$ be a reduced multi-curve with Dehn-Thurston coordinates
\beq\lb{2s}
n_i(\Omega)=\begin{cases} 0 & \text{for}\  i=1,...,u,\\
2 & \text{for}\ i=u+1,..., d+s,\\
\end{cases}
\eeq
with all of twist coordinates zero. The distribution of the intersection numbers $n_i=2$ for $u< i\leq d$ among $n_b$ for $b\in B$ in the $i$-th boundary component can be chosen arbitrarily here.
Note that $\Omega$ exists. (It is the boundary of a regular neighborhood of the graph composed of all base intervals in annuli $\mathbb A_i$ and of all $Y$ graphs in pairs of pants $\mathbb P_{\{i,j,k\}}$.)
$\Omega\in \cal H_{0,...,0,2,...,2}$ (where the first $u$ indexes are zero).

Assuming that $c$ involves a monomial as above, we are going to show that $\Omega\cdot c-c\cdot \Omega\ne 0.$
For any monomial 
\beq \lb{e_c-mon}
m=r\cdot \gamma_1^{t_1}\cdot ...\cdot \gamma_{d+s}^{t_{d+s}},
\eeq
in $c$ 
the values of $m\cdot \Omega$ and of $\Omega\cdot m$ in $\cal{GRS}(F,B)$ can be computed by considering a contribution of each annulus separately by 
Remark \ref{r-gen-con}.
The left and right multiplications of $\Omega$ by $\gamma_1^{t_1}\cdot ...\cdot \gamma_{u}^{t_u}$ change the first $u$
components of $\Omega$ to $t_1,...,t_u$.
The left and right multiplications of $\Omega$ by $\gamma_i^{t_i}$ for $i>u$ can be computed by (\ref{e-annulus-l}) and (\ref{e-annulus-r}). They will be linear combinations of reduced multi-curves with all twist coordinates zero, except possibly for the $i$-th one. Putting the effects of multiplications by $\gamma_i^{t_i}$ for $i=1,...,d+s$ together, we conclude that
$\Omega\cdot m-m\cdot \Omega$ is an element of $\cal H_{0,...,0,2,...,2}$ corresponding to
$$r\prod_{i=1}^{u} x_i^{t_i}\cdot \left( \prod_{i=u+1}^{d+s} \sum_{j=0}^{t_i} {t_i \choose j} A^{2t_i-4j} x_i^{t_i-2j}-
\prod_{i=u+1}^{d+s} \sum_{j=0}^{t_i} {t_i \choose j} A^{-2t_i+4j} x_i^{t_i-2j}\right).$$

Now assume that (\ref{e_c-mon}) is the leading monomial, $LM(c)$, 
of the polynomial $c$ with respect to the lexicographic ordering of the variables gamma, such that
$\gamma_1< .... <\gamma_{d+s}.$ Then by the above formula,
$$LM(\Omega\cdot c- c \cdot \Omega)=r\prod_{i=1}^{u} x_i^{t_i} \left(\prod_{i=u+1}^{d+s}
A^{2t_i} x_i^{t_i} - \prod_{i=u+1}^{d+s}
A^{-2t_i} x_i^{t_i}\right)=$$
$$r\left(A^{2\sum_{i=u+1}^{d+s} t_i}  - A^{-2\sum_{i=u+1}^{d+s} t_i}\right) \prod_{i=1}^{d+s} x_i^{t_i}.$$
By our assumption of $t_i>0$ for some $i>u,$ this expression is non-zero, since $A^{4\sum_{i=u+1}^{d+s} t_i}  -1$ is not a zero divisor in $R$. 
Therefore, $c$ cannot be central. That completes the argument by contradiction.

%
\section{Proof of Theorem \ref{main-div-gr}}
\lb{s_proof-div}
%

Consider a total order $\prec$ on $\cal{RMC}(F,B)$ such that $C \prec C'$ if (a) $C$ is of lower weight than $C'$ or (b) $C$ and $C'$ have equal weight but
$C$ is less than $C'$ in the lexicographic order on $\cal{RMC}(F,B)$ identified by Corollary \ref{D-T-B-cor} with a subset of $\Z_{\geq 0}^{|B|+s}\times \Z^{d+s}$ (through Dehn-Thurston coordinates).

For any $x=\sum_{C\in \cal{RMC}(F,B)} r_C\, C \in \cal{GRS}(F,B)$, let $LM(x)$ be the largest summand of $x$ with respect of that ordering, that is 
$$LM(x)=r_{C_0}C_0,$$
where $C_0$ is the largest $C$ with non-zero coefficient $r_C.$

\begin{figure}
$\parbox{1.3in}{\includegraphics[height=1.3in]{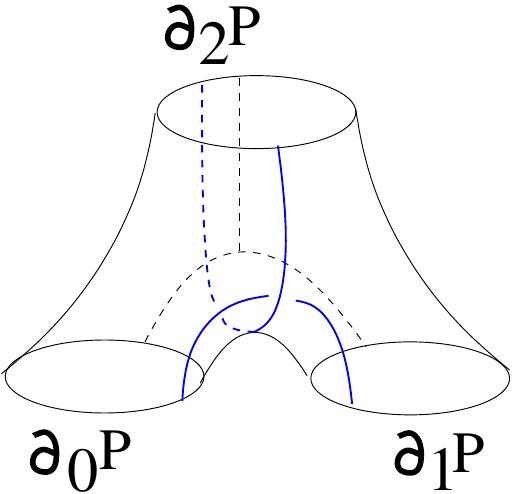}}= A
\parbox{1.3in}{\includegraphics[height=1.3in]{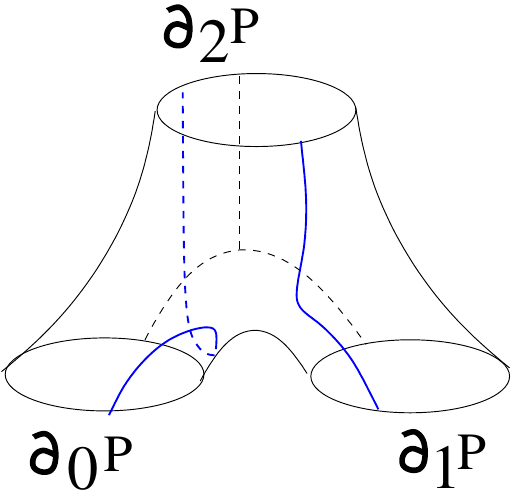}} +A^{-1}
\parbox{1.3in}{\includegraphics[height=1.3in]{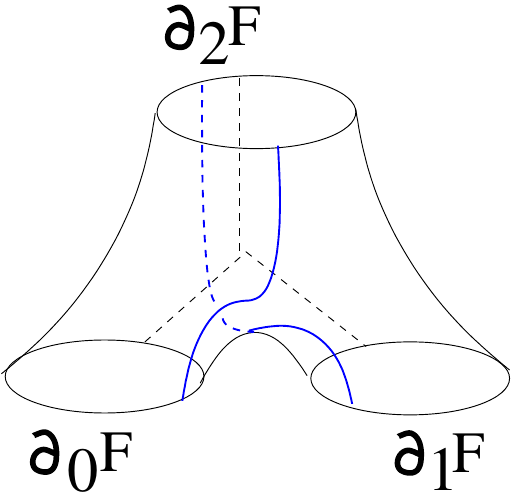}}$
$A\parbox{1.3in}{\includegraphics[height=1.3in]{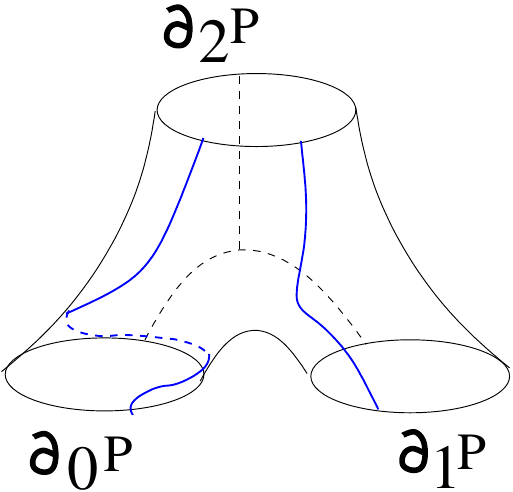}} +A^{-1}\parbox{1.3in}{\includegraphics[height=1.3in]{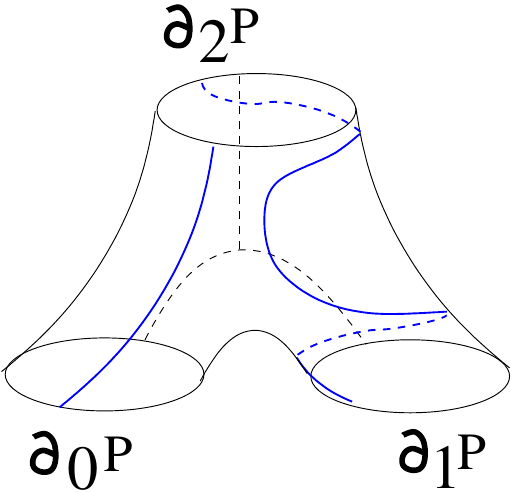}}$
\caption{$DT(1,1,0,0,0,0)\cdot DT(0,0,2,0,0,0,0)$ and its $A$ resolution, $DT(1,1,2,1,0,0)$, and
$A^{-1}$ resolution, $DT(1,1,2,0,-1,1).$}
\lb{eg-trouble}
\end{figure}

From now one the following notation will be useful: For any $x,x'$ in some $R$-algebra we will write $x \dot{=} x'$ if $x= A^e x'$ for some $e\in \Z.$

If for every $x,y\in \cal{GRS}(F,B)$, 
\beq\lb{e-hope}
LM(x\cdot y)\dot{=} LM(x)\cdot LM(y)
\eeq
was to hold in $\cal{GRS}(F,B)$, that would imply that $\cal{GRS}(F,B)$ is a domain for every integral domain $R$. That is not the case however, as example of Fig. \ref{eg-trouble} shows. Denote by $D(n_0,n_1,n_2,t_0,t_1,t_2)$ the multi-curve in $(\mathbb P, \p \mathbb P\cap Y)$ with these
Dehn-Thurston coordinates. Then\vspace*{.1in}\\
$DT(1,1,0,0,0,0)\cdot DT(0,0,2,0,0,0,0)=$\hfill
$$A\cdot DT(1,1,2,1,0,0)+A^{-1}\cdot DT(1,1,2, 0,-1,1)$$
in $\cal{GRS}(\mathbb P,\p \mathbb P\cap Y)$.\footnote{The asymmetry of torsion may look surprising, but recall that we always draw surfaces so that they are oriented counterclockwise ``up front''. Hence, taking mirror image does not create inconsistency here.}
The above multiplication does not satisfy (\ref{e-hope}).
There is, however, a nice subalgebra of $\cal{GRS}(F,B)$ 
in which (\ref{e-hope}) holds, defined as follows:

As before, let $n_i(C)=\sum n_b(C)$, for $i=1,...,d$, where the sum is over $b\in B$ lying 
on the $i$-th boundary component of $\p F.$
We say that $C\in \cal{RMC}(F,B)$ is {\em triangular} if 
\begin{itemize}
\item for every pair of pants
$\mathbb P_{\{i_0,i_1,i_2\}}\in \cal P,$ $n_{i_0}(C), n_{i_1}(C), n_{i_2}(C)$ satisfy the triangle inequalities (\ref{e_triang-ineq}). 
\item $t_i=0$ if $n_i=0.$
\end{itemize}
Denote the set of triangular reduced multi-curves by $\cal{RMC}^\triangle(F,B)$ and the subspace of $\cal{GRS}(F,B)$ spanned by them by $\cal{GRS}^\triangle(F,B).$ 

\bpro\lb{pro-triang} 
For every $C,C' \in \cal{RMC}^\triangle(F,B)$ 
$$C\cdot C' \dot{=} C''$$ 
in $\cal{GRS}(F,B),$
where $C''$ is a triangular multi-curve with Dehn-Thurston coordinates
being the sum of those of $C$ and $C',$ 
$\nu''=\nu +\nu'$, $n_i''=n_i+n_i'$, for $i=d+1,...,d+s,$ 
and $t_i''=t_i+t_i'$, for $i=1,...,d+s.$
\epro

The proof of Proposition \ref{pro-triang} below holds for any $B$, but it is worth noting that it is enough to prove this result for $B'$ instead (which is slightly simpler because $B'$ has at most one base point per component of $\p F$). Indeed, there is a natural ``forgetful'' $R$-module epimorphism: 
$$\Psi: \cal{GRS}(F,B)\to \cal{GRS}(F,B').$$ 
Note that for every $C,C'\in \cal{RMC}(F,B)$,
\beq\lb{Psi-mult}
\Psi(C\cdot C') \dot{=} \Psi(C) \cdot \Psi(C'),
\eeq
because the diagrams of $\Psi(C\cdot C')$ and of $\Psi(C) \cdot \Psi(C')$ differ by some possible crossing changes of arcs near $\p F$ only, and only one resolution per each such crossing in non-zero in the reduced skein module.
Suppose now that Proposition \ref{pro-triang} holds for $(F,B')$.
Observe that $\cal{GRS}(F,B)$ is multi-graded by
$\nu\in \Z_{\geq 0}^B$, and that $\Psi$  is $1$-$1$ on each homogenous component of that grading. 
Therefore, if $C,C' \in \cal{RMC}^\triangle(F,B)$ and (\ref{Psi-mult}) holds then since  $C\cdot C'$ and $C''$ belong to the same homogeneous component of $\cal{GRS}(F,B)$, we have 
$C\cdot C' \dot{=} C''$ (by the above observation).\vspace*{.1in}

\noindent{\it Proof of Proposition \ref{pro-triang}:}
%
%
$$(C\cdot C')\cap \mathbb A_i \dot{=} (C''\cap \mathbb A_i),$$
in $\cal{GRS}_{n_1+n_2}(\mathbb A_i)$ for every $i=1,...,d+s$, by Lemma \ref{lem-elem-tang-prod}. (If $\mathbb A_i$ is a boundary annulus of $\p F$, with marked points $b_1,...,b_k$, $k>1$, on the corresponding component of $\p F$ then one needs to consider $\cal{GRS}_{n_1+n_2}(\mathbb A_i,\{*\}\times \{1\}\cup \{b_1,...,b_k\}\times \{-1\})$ instead, but the same formula holds.)

Also note that $(C\cdot C')\cap \mathbb P_{\{i,j,k\}}$ can be transformed into 
$$T_{i,j,k}(n_i(C)+n_i(C'),n_j(C)+n_j(C'),n_k(C)+n_k(C'))$$ 
by rearranging endpoints of
$C\cap \mathbb P_{\{i,j,k\}}$ and of $C'\cap \mathbb P_{\{i,j,k\}}$, i.e. by undoing crossings between them.
That undoing of crossing (or, equivalently, adding opposite braiding of endpoints) changes the value of a diagram in $\cal{GRS}_{n_1+n_2}(A_i)$ by a power of $A.$
Since $C\cdot C'$ and $C''$ differ by powers of $A$ only in every annulus and every pair of pants,
the statement now follows from Lemma \ref{l_conc} and Remark \ref{r-gen-con}.
\qed


\bcor $\cal{GRS}^\triangle(F,B)$ is an $R$-subalgebra of $\cal{GRS}(F,B)$.
\ecor

\bpr
If $n_{i_0}(C), n_{i_1}(C), n_{i_2}(C)$ and $n_{i_0}(C'), n_{i_1}(C'), n_{i_2}(C')$ satisfy the triangle inequalities then 
$n_{i_0}(C)+n_{i_0}(C'), n_{i_1}(C)+n_{i_1}(C'), n_{i_2}(C)+n_{i_2}(C')$ do to. Therefore, the statement follows from the proposition above.
\epr

\bcor \lb{triangle-no-zero}
$\cal{GRS}^\triangle(F,B)$ is a domain. 
\ecor 

\bpr
By the above proposition, for any non-zero $x,y\in \cal{GRS}^\triangle(F,B)$, 
$$LM(x\cdot y) \dot{=} LM(x)\cdot LM(y).$$ 
Furthermore, by the same proposition, $LM(x)\cdot LM(y)$ is a non-zero scalar multiple of a reduced multi-curve in $\cal{GRS}^\triangle(F,B)$ and, hence, it is non-zero in $\cal{GRS}^\triangle(F,B)$.
\epr

Now we are going to deduce Theorem \ref{main-div-gr} from Corollary \ref{triangle-no-zero} and the following statement. 
Let $\Omega\in \cal{RMC}(F,B)$ be as in the previous section, i.e. have intersection coordinates (\ref{2s}) and all twist coordinates zero.
(If $B$ has more than a single point on some component of $\p F$ then such $\Omega$ is not unique.)

Let $r_{\Omega}, l_{\Omega}: \cal{GRS}(F,B) \to \cal{GRS}(F,B)$ be the right and left multiplication maps, $r_{\Omega}(x)=x\cdot \Omega,\ l_{\Omega}(x)=\Omega \cdot x$.

\bpro \lb{T0-to-triang} 
(1) For every $x\in \cal{GRS}(F,B)$, 
$$r_{\Omega}^n(x),l_{\Omega}^n(x)\in \cal{GRS}^\triangle(F,B),$$ 
for sufficiently large $n$.\\
(2) The maps $r_{\Omega}, l_{\Omega}: \cal{GRS}(F,B) \to \cal{GRS}(F,B)$
are $1$-$1$.
\epro 


\noindent {\bf Proof of Theorem \ref{main-div-gr}:}
Let $x,y\in \cal{GRS}(F,B),$ $x,y\ne 0.$ By the above proposition, $\Omega^n x$ and $y \Omega^n$ are non-zero elements of $\cal{GRS}^\triangle(F,B)$ for sufficiently large $n$ and, hence, $\Omega^n x y \Omega^n\ne 0$ in $\cal{GRS}^\triangle(F,B)$ by Corollary \ref{triangle-no-zero}. That implies in particular that $xy\ne 0$ in
$\cal{GRS}(F,B)$.
\qed\vspace*{.1in}


For the purpose of proving Proposition \ref{T0-to-triang} we will need the following:

\blem \lb{lemma-pants-tc} 
Assume that $n_0,n_1,n_2\geq 0,$ $n_0+n_1+n_2$ is even, and
$n_i>n_{i+1}+n_{i+2}$, for some $i=0,1,2,$ where the indexes are mod $3$.
Then  
$$T(n_0,n_1,n_2)\cdot T(2,2,2) \dot{=} (A^{n_i'}\cdot  D_1+ A^{-n_i'}\cdot  D_2),$$
in $\cal{RS}_{n_0+2,n_1+2,n_2+2}(\mathbb P)$, where $n_i'=\frac{n_i-n_{i+1}-n_{i+2}}{2}$ and we identity $\mathbb P$ with $\mathbb P\cup \mathbb A_0\cup \mathbb A_1\cup \mathbb A_2$ (adjacent annuli to $\p_0 \mathbb P,$ $\p_1 \mathbb P$, and $\p_2 \mathbb P$, respectively) and
\begin{itemize} 
\item $D_1$ is $T(n_0+2,n_1+2,n_2+2)$ 
composed with $T(n_{i+1}+2,1)$ in $\mathbb A_{i+1}$. 
\item $D_2$ is $T(n_0+2,n_1+2,n_2+2)$ composed with $T(n_i+2, 1)$ in  $\mathbb A_i$ and with 
$T(n_{i+2}+2, -1)$ in  $\mathbb A_{i+2}$.
\end{itemize}
\elem 

\begin{figure}
\includegraphics[height=1.4in]{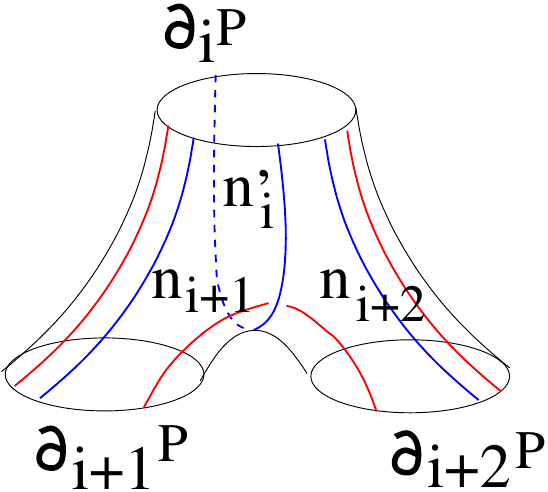}
\caption{Product of $T(n_0,n_1,n_2)$ (in blue) with $T(2,2,2)$ (in red), for $n_i>n_{i+1}+n_{i+2}$, with its endpoints rearranged to minimize the crossing number. For readability of the picture, we did not depict the $Y$ graph in the back of $\mathbb P$. 
}
\lb{pants-tc}
\end{figure} 

Notice that the above lemma is a generalization of the multiplication in Fig. \ref{eg-trouble}.

\bpr  

Since crossing changes between link endpoints in $\cal{RS}_{n_0+2,n_1+2,n_2+2}(\mathbb P)$ change the value of the link by a multiplicative factor of power of $A$, we can consider Fig. \ref{pants-tc} instead of $T(n_0,n_1,n_2)\cdot T(2,2,2)$ for the purpose of the proof. Its $n_i'$ crossings are enclosed in a subdiagram in Fig. \ref{subdiagram}(a).
\begin{figure}
\includegraphics[height=0.7in]{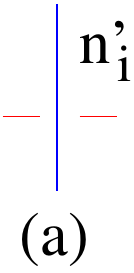}\hspace*{.2in}
\includegraphics[height=0.7in]{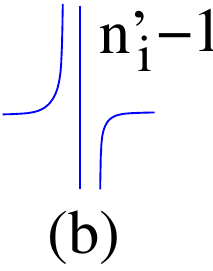}\hspace*{.2in}\ 
\includegraphics[height=0.7in]{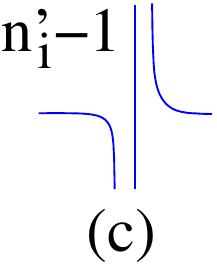}
\caption{}
\lb{subdiagram}
\end{figure} 
All resolutions of these crossings result in cups or caps in the subdiagram and, hence, vanish in $\cal{RS}_{n_0+2,n_1+2,n_2+2}(\mathbb P)$, except for the one with all resolutions of type $A$ and the one with all resolutions of type $A^{-1}$ which yield the subdiagrams shown in Fig. \ref{subdiagram}(b) and (c) respectively and the diagrams $D_1$ and $D_2$ in Fig. \ref{D1D2}.  
\epr

\begin{figure}
\parbox{1.4in}{\includegraphics[height=1.2in]{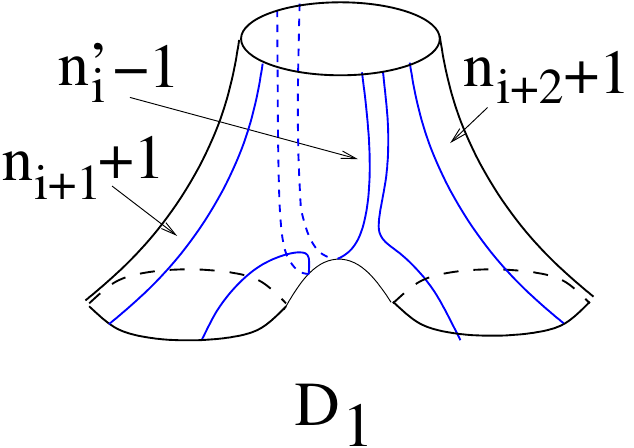}} \hspace*{.2in} $=$ \hspace*{.1in}
\parbox{1.4in}{\includegraphics[height=1in]{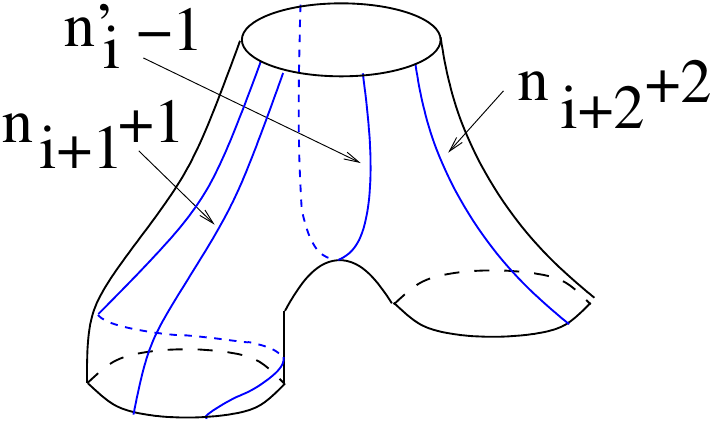}}\vspace{.2in}

\parbox{1.6in}{\includegraphics[height=1.1in]{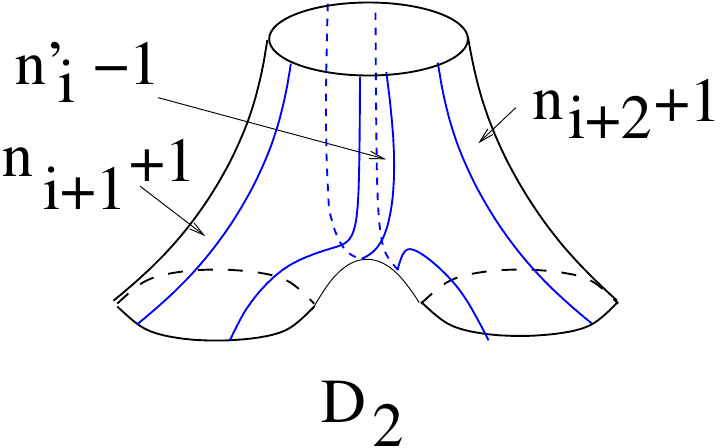}} $=$ \hspace*{.1in}
\parbox{1.6in}{\includegraphics[height=1in]{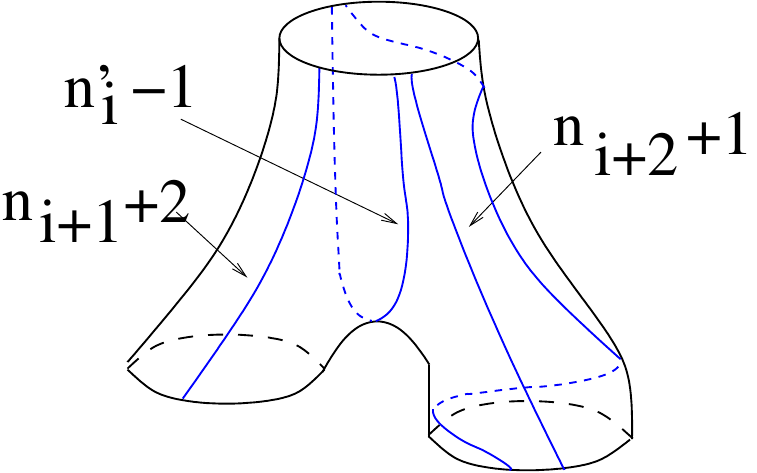}}
\caption{Diagrams $D_1$ and $D_2$.}
\lb{D1D2}
\end{figure}

\noindent{\bf Proof of Proposition \ref{T0-to-triang}:}

(1) $r_{\Omega}^n$ and $l_{\Omega}^n$ increase all intersection coordinates $n_i$, for $i>u,$ of any reduced multi-curve $C$ in $F$ by $2n$. Now the statement follows from the fact that there are finitely many pairs of pants in $F$ and for each pair of pants $\mathbb P_{\{i,j,k\}}$
$$n_i(C)+2n,n_j(C)+2n,n_k(C)+2n$$ 
satisfy triangle inequalities for sufficiently large $n$.

(2) The proof below works for arbitrary $B$. Alternatively, one can reduce the statement to that for $(F,B')$ by the same method as in the remark following the statement of Proposition \ref{pro-triang}.

As in the previous section, let $\cal H_{\nu, n_{d+1},...,n_{d+s}}$ be the free $R$-module with a basis given by reduced multi-curves $C\in \cal{RMC}(F,B)$ with intersection coordinates $\nu, n_{d+1},...,n_{d+s}.$ Recall that such $\cal H$-spaces define an algebra grading on $\cal{GRS}(F,B)$.  Since 
$l_{\Omega}, r_{\Omega}$ map different $\cal H$-spaces to different $\cal H$-spaces, it is enough to show that $l_{\Omega}, r_{\Omega}$ are $1$-$1$ when restricted to $\cal H$-spaces.

As in the previous section, assume that $\gamma_1,...,\gamma_u$ are the unmarked components of $\p F$. 
The functions $r_{\Omega}$ and $l_{\Omega}$ map $\cal H_{\nu, n_{d+1},...,n_{d+s}}$ to\\$\cal H_{\nu+\nu_\Omega, n_{d+1}+2,...,n_{d+s}+2}$ where 
the indices of the $H$-space on the right come from the sum of the intersection coordinates $(\nu, n_{d+1},...,n_{d+s})$ with the intersection coordinates of $\Omega,$
$(\nu_\Omega, 2,....,2)$. Here, $\nu_\Omega(b)=2$ for one $b\in B$ on each marked boundary component of $\p F$ and $\nu_\Omega(b')=0$ for all other $b'\in B.$
Assume that $(\nu, n_{d+1},...,n_{d+s})$ is realizable.
By identifying $\cal H_{\nu, n_{d+1},...,n_{d+s}}$ and $\cal H_{\nu+\nu_\Omega, n_{d+1}+2,...,n_{d+s}+2}$
with $R[x_1,...,x_u, x_{i+1}^{\pm 1},..., x_{d+s}^{\pm 1}]$ as before, 
we consider $r_{\Omega}, l_{\Omega}$ as functions 
$$r_{\Omega}, l_{\Omega}: R[x_1,...,x_u, x_{u+1}^{\pm 1},..., 
x_{d+s}^{\pm 1}]\to R[x_1,...,x_u, x_{u+1}^{\pm 1},..., x_{d+s}^{\pm 1}].$$

Consider a pair of pants $\mathbb P_{i,j,k}$ such that $n_i,n_j,n_k$ fail some triangle inequality. Without loss of generality, we can assume that $n_i>n_j+n_k.$ 
By its definition, $\mathbb P_{i,j,k}$ is bounded by $\gamma_{i,\ve_1},\gamma_{j,\ve_2},\gamma_{k,\ve_3}$ for some $\ve_1,\ve_2,\ve_3\in \{\pm\}.$ 
Assume that if $\mathbb P_{i,j,k}$ is drawn such that $Y$ graph is in the back and the front is oriented counterclockwise then $\gamma_{i,\ve_1},\gamma_{j,\ve_2},\gamma_{k,\ve_3}$ appear in counterclockwise direction.
Let
$$p_{i,j,k}=A^{n_i'} x_j+A^{-n_i'}x_ix_k^{-1},$$
where $n_i'=\frac{n_i-n_j-n_k}{2}$.

\blem \lb{lem-rT0} 
For every reduced multi-curve $C$ in $\cal{GRS}(F,B)$ with intersection coordinates $\nu, n_{d+1},...,n_{d+s},$
$$r_{\Omega}(C) \dot{=} \prod_{i=u+1}^{d+s} y_i \cdot \prod p_{i,j,k} \in R[x_1,...,x_u, x_{u+1}^{\pm 1},..., x_{d+s}^{\pm 1}]\simeq \cal H_{\nu+\nu_\Omega, n_d+2,...,n_{d+s}+2},$$
where $y_i=x_i^{t_i(C)}$
if $n_i\ne 0$ and 
$$y_i= (A^2 x_i+A^{-2} x_i^{-1})^{t_i(C)}$$ 
for $i$ such that $n_i=0.$
The product $\prod p_{i,j,k}$ is over all pairs of pants $\cal P_{\{i,j,k\}}$ such that $n_i,n_j,n_k$ fail to satisfy triangle inequalities. (The $i$ in $p_{i,j,k}$ is unrelated to the $i$ in $y_i.$)
\elem

\bpr For any $C\in \cal{RMC}(F,B)$,
$C\cdot \Omega$ is given by the concatenation $\Phi$' (of Lemma \ref{l_conc} and of Remark \ref{r-gen-con}) of the pieces $(C\cdot \Omega)\cap \mathbb A_i$, for $i=1,...,d+s$, with the pieces $(C\cdot \Omega)\cap \mathbb P$ for all pairs of pants $\mathbb P\in \cal P$. 

By Lemma \ref{lem-elem-tang-prod} and by (\ref{e-annulus-l}), (\ref{e-annulus-r}),
$$(C\cdot \Omega)\cap \mathbb A_i \dot{=}
\begin{cases} T(n_i+2,t_i(C)) & \text{for}\ n_i\ne 0,\\
\sum_{l=1}^{t_i} {t_i(C) \choose l} A^{2t_i(C)-4l} T(2,t_i(C)-2l)
& \text{for}\ n_i= 0,\\
\end{cases}
$$
for $i=u+1,...,d+s.$ 
Given the identifications
\beq\lb{e-H-spaces}
\cal H_{\nu, n_d,...,n_{d+s}}\simeq  R[x_1,...,x_u, x_{u+1}^{\pm 1},..., x_{d+s}^{\pm 1}]\simeq \cal H_{\nu+\nu_\Omega, n_d+2,...,n_{d+s}+2},
\eeq
the above multiplication operation corresponds either to identity map or
to multiplication by
$$\sum_{l=0}^{t_i(C)} {t_i(C) \choose l} A^{2 t_i(C)-4l} x_i^{t_i(C)-2l}=
(A^2 x_i+A^{-2} x_i^{-1})^{t_i(C)}.$$ 

Let us consider the effect of the multiplication of $C$ by $\Omega$ in pairs of pants now. By assuming $C$ is in a canonical position, $C\cap \mathbb P_{\{i,j,k\}}= T(n_i,n_j,n_k)$ for some $n_i,n_j,n_k.$ If they satisfy triangle inequalities, then 
$$(C\cdot \Omega)\cap \mathbb P_{\{i,j,k\}} \dot{=} T(n_i+2,n_j+2,n_k+2)$$ 
in $\cal{RS}_{n_i+2,n_j+2,n_k+2}(\mathbb P)$ by Proposition \ref{pro-triang}. 

If $n_i,n_j,n_k$ fail some triangle inequality, then the value of $(C\cdot \Omega)\cap \mathbb P$ in $\cal{RS}_{n_i+2,n_j+2,n_k+2}(\mathbb P)$ is given by Lemma \ref{lemma-pants-tc}.

Since twist coordinates are additive under composition of canonical link diagrams in annuli,  
the $i$-th twist coordinates of the summands of $C\cdot \Omega$ are the sums of twists coming from the multiplication $C\cdot \Omega$ in $\mathbb A_i$ and of the possible correction twists from the adjacent pairs of pants, as in Lemma \ref{lemma-pants-tc}.
Therefore, through identification (\ref{e-H-spaces}) the twist coordinates are given by the formula of the statement.
\epr 

Going back to the proof of Proposition \ref{T0-to-triang}, observe that the term $\prod p_{i,j,k}$ in the formula of Lemma \ref{lem-rT0} does not depend on 
the twist coordinates of $C$. In other words, this product depends on the values of $n_1,...,n_{d+s}$ only.
Therefore, it is enough to prove that 
$$\bar r_{\Omega}=r_{\Omega}/\prod p_{i,j,k}: R[x_1,...,x_u, x_{u+1}^{\pm 1},..., x_{d+s}^{\pm 1}]\to 
R[x_1,...,x_u, x_{u+1}^{\pm 1},..., x_{d+s}^{\pm 1}]$$
is $1$-$1$.
By Lemma \ref{lem-rT0}, 
$$\bar r_{\Omega}(x_1^{t_1}\cdot ...\cdot x_{d+s}^{t_{d+s}}) \dot{=} \cdot \prod_{i=1}^{d+s} y_i ,$$
where $y_i$ are as above. 
Therefore, given a lexicographic ordering of monomials in $x_1,...,x_{d+s}$, 
the leading monomial of $\bar r_{\Omega}(x_1^{t_1}\cdot ...\cdot x_{d+s}^{t_{d+s}})$ with respect to that ordering is $x_1^{t_1}\cdot ...\cdot x_{d+s}^{t_{d+s}}$ up to a power of $A.$ That means that for any polynomial $w$, the leading terms of $w$ and of $\bar r_{\Omega}(w)$ 
coincide up to a power of $A.$ Therefore,  $\bar r_{\Omega}$ is $1$-$1$ and, consequently, $r_{\Omega}$ is $1$-$1$ as well.

The proof of $l_{\Omega}$ being $1$-$1$ is analogous. (Note that $\Omega\cdot C$ becomes $C\cdot \Omega$ by substituting $A$ by $A^{-1}$.)
That completes the proof of Proposition \ref{T0-to-triang} and of Theorem \ref{main-div-gr}.
\qed

%

\end{document}